\documentclass{article}

\usepackage[]{arXiv}
\usepackage[table,xcdraw]{xcolor} 
\usepackage{scalerel}
\usepackage{tikz}
\usetikzlibrary{svg.path}

\definecolor{orcidlogocol}{HTML}{A6CE39}
\tikzset{
  orcidlogo/.pic={
    \fill[orcidlogocol] svg{M256,128c0,70.7-57.3,128-128,128C57.3,256,0,198.7,0,128C0,57.3,57.3,0,128,0C198.7,0,256,57.3,256,128z};
    \fill[white] svg{M86.3,186.2H70.9V79.1h15.4v48.4V186.2z}
                 svg{M108.9,79.1h41.6c39.6,0,57,28.3,57,53.6c0,27.5-21.5,53.6-56.8,53.6h-41.8V79.1z M124.3,172.4h24.5c34.9,0,42.9-26.5,42.9-39.7c0-21.5-13.7-39.7-43.7-39.7h-23.7V172.4z}
                 svg{M88.7,56.8c0,5.5-4.5,10.1-10.1,10.1c-5.6,0-10.1-4.6-10.1-10.1c0-5.6,4.5-10.1,10.1-10.1C84.2,46.7,88.7,51.3,88.7,56.8z};
  }
}

\newcommand\orcidicon[1]{\href{https://orcid.org/#1}{\mbox{\scalerel*{
\begin{tikzpicture}[yscale=-1,transform shape]
\pic{orcidlogo};
\end{tikzpicture}
}{|}}}}

\usepackage[numbers,sort&compress]{natbib}
\bibpunct[, ]{[}{]}{,}{n}{,}{,}

\usepackage{graphicx}
\usepackage[utf8]{inputenc}
\usepackage{amsfonts}
\usepackage{tabularx,booktabs}

\usepackage{amsmath}
\usepackage{mathtools}
\usepackage{bbm}
\usepackage{makecell}

\newcommand{\argmin}{\operatornamewithlimits{argmin}}

\usepackage{amsmath,amssymb,amsfonts}
\usepackage{textcomp,setspace}

\usepackage{float}
\usepackage[utf8]{inputenc} 
\usepackage[T1]{fontenc}    
\usepackage{nth}
\usepackage{algorithm,algpseudocode}      
\usepackage{url}            
\usepackage{booktabs}       
\usepackage{amsfonts}       
\usepackage{nicefrac}       
\usepackage{microtype}      
\usepackage{array}
\usepackage{amsbsy} 
\usepackage{amsmath}
\usepackage{ulem}
\usepackage{blkarray}
\usepackage{amsmath}
\newtheorem{theorem}{Theorem}

\def\BibTeX{{\rm B\kern-.05em{\sc i\kern-.025em b}\kern-.08em
		T\kern-.1667em\lower.7ex\hbox{E}\kern-.125emX}}

\makeatletter
\def\BState{\State\hskip-\ALG@thistlm}
\makeatother

\newcommand{\st}{{\rm s.t.}}
\newcommand{\cG}{\mathcal{G}}
\newcommand{\cV}{\mathcal{V}}
\newcommand{\cE}{\mathcal{E}}

\newcommand{\cN}{\mathcal{N}}
\newcommand{\cI}{\mathcal{I}}
\newcommand{\cR}{\mathcal{R}}

\newcommand{\br}{\mathbf{r}}
\newcommand{\bR}{\mathbf{R}}
\newcommand{\bT}{\mathbf{T}}
\newcommand{\bv}{\mathbf{v}}
\newcommand{\bs}{\mathbf{s}}
\newcommand{\bS}{\mathbf{S}}
\newcommand{\ba}{\mathbf{a}}
\newcommand{\bA}{\mathbf{A}}
\newcommand{\bc}{\mathbf{c}}
\newcommand{\bq}{\mathbf{q}}
\newcommand{\bD}{\mathbf{D}}
\newcommand{\bP}{\mathbf{P}}

\newcommand{\bu}{\mathbf{u}}
\newcommand{\bH}{\mathbf{H}}
\newcommand{\bW}{\mathbf{W}}

\newcommand{\bI}{\mathbf{I}}
\newcommand{\bz}{\mathbf{z}}
\newcommand{\bga}{\boldsymbol{\gamma}}

\newcommand{\bla}{\boldsymbol{\lambda}}
\newcommand{\bLa}{\boldsymbol{\Lambda}}

\newcommand{\bbR}{\mathbb{R}}

\usepackage{eso-pic}
\usepackage{multirow}

\usepackage{times}
\usepackage{fancyhdr,graphicx,amsmath,amssymb}
\include{pythonlisting}
\usepackage{etoolbox, lineno}
\usepackage{lipsum}
\usepackage{url}

\usepackage{breakurl}
\usepackage[breaklinks]{hyperref}
\begin{document}

\title{Congestion Reduction via Personalized Incentives}

\author{Ali~Ghafelebashi \orcidicon{0000-0001-8339-7960},
        Meisam~Razaviyayn \orcidicon{0000-0003-4342-6661},
        and~Maged~Dessouky \orcidicon{0000-0002-9630-6201}
\thanks{This study was funded by a grant from the National Center for Sustainable Transportation (NCST), supported by the U.S. Department of Transportation’s University Transportation Centers Program. The contents of this project reflect the views of the authors, who are responsible for the facts and the accuracy of the information presented herein. This document is disseminated in the interest of information exchange. The U.S. Government and the State of California assumes no liability for the contents or use thereof. Nor does the content necessarily reflect the official views or policies of the U.S. Government and the State of California. This paper does not constitute a standard, specification, or regulation. This paper does not constitute an endorsement by the California Department of Transportation (Caltrans) of any product described herein.}
}

\maketitle

\begin{abstract}
With rapid population growth and urban development, traffic congestion has become an inescapable issue, especially in large cities. Many congestion reduction strategies have been proposed in the past, ranging from roadway extension to transportation demand management. In particular, congestion pricing schemes have been used as negative reinforcements for traffic control. In this project, we study an alternative approach of offering positive incentives to drivers to take different routes. More specifically, we propose an algorithm to reduce traffic congestion and improve routing efficiency via offering personalized incentives to drivers.  We exploit the wide-accessibility of smart devices to communicate with drivers and develop an incentive offering mechanism using individuals’ preferences and aggregate traffic information. The incentives are offered after solving a large-scale optimization problem in order to minimize the total travel time (or minimize any cost function of the network such as total Carbon emission). Since this massive size optimization problem needs to be solved continually in the network, we developed a distributed computational approach. The proposed distributed algorithm is guaranteed to converge under a mild set of assumptions that are verified with real data. We evaluated the performance of our algorithm using traffic data from the Los Angeles area. Our experiments show congestion reduction of up to 11\% in arterial roads and highways.
\end{abstract}


\section{Introduction}
\label{sec:Introduction}
Today, traffic congestion is one of the most prevalent issues in large metropolitan areas, resulting in lowered quality of life for residents and economic losses. According to INRIX the United States' economy suffered \$88 billion in losses and the drivers lost 99 hours due to traffic congestion in 2019~\cite{INRIX}. In addition to direct economic losses, traffic congestion can worsen air quality and adversely affect health conditions. According to the Transportation Research Board, vehicle emission is the main cause for air pollution~\cite{national2002congestion, health2010traffic, zhang2013air, hennessy1999traffic, selzer1974life}. 

\vspace{0.2cm}

In this paper, the focus is on changing drivers' behavior by offering incentives in order to reduce traffic congestion. Perhaps one of the closest strategies to our solution is pricing mechanisms in the literature. Road pricing policies, such as assigning a fee or tax for driving on a highway/road, have been widely studied in theory and practice~\cite{pigou1920economics, knight1924some, RePEc:elg:eebook:4192, van2009behavioural, bouchelaghem2018reliable, zhang2013self, kachroo2016optimal, farokhi2014study, groot2014toward, cao2020improving}. In this approach, congestion reduction is expected as a result of discouraging people to use congested roads. Such pricing mechanisms could be dependent on different factors such as time~\cite{zheng2016time}, distance~\cite{daganzo2015distance}, or vehicle characteristics~\cite{zhang2019impact}. 
While pricing is a promising approach from a market point of view, issues such as equity barriers complicate the implementation of congestion pricing/taxation schemes~\cite{knockaert2012spitsmijden, levinson2010equity, martens2012justice, raux2004acceptability, ieromonachou2006evaluation, hensher2014type}. In addition, complexities and uncertainties in designing pricing mechanisms have prevented policymakers from implementing advanced congestion pricing schemes~\cite{gu2018congestion}. 
Tradable credits (TCs) or tradable mobility permits (TPMs) are another token-based pricing mechanism~\cite{verhoef1997tradeable, viegas2001making, raux2004use, fan2013tradable, thogersen2008breaking, wang2014models}. The theoretical advantages of such tradable credits have been shown in~\cite{tsekeris2009design, nie2015new, wu2012design, grant2014role}. While such cap-and-trade programs have been implemented in some economic sectors, such as airport slot allocation~\cite{fukui2010empirical}, it has not been implemented for individual-level personal travels and daily commutes~\cite{dogterom2017tradable} due to the design complexities of such token markets~\cite{azevedo2018tripod}. 
\vspace{0.2cm}

Lately, researchers have paid more attention to positive incentive policies. Based on the psychological theory of reactance, rewarding desirable behavior could work better than penalizing undesirable behavior~\cite{Brehm1966-BREATO-2}. Moreover, rewarding is a more popular policy than a taxation approach~\cite{knockaert2007experimental}. While the effectiveness of rewarding in changing the individual's behavior have been shown in \cite{kreps1997intrinsic} and \cite{berridge2001reward}, there is a limited number of studies on the effectiveness of rewarding policies in the transportation area. Among these studies, the INSTANT project~\cite{merugu2009incentive} and the CAPRI project~\cite{yue2015reducing} have provided positive incentives to motivate commuters to avoid peak time for their travel. They have shown the effectiveness of a rewarding policy in congestion reduction. Another series of studies has been done in the Netherlands on the effectiveness of monetary rewards to avoid rush hour driving~\cite{bliemer2009rewarding, ben2011behaviour, ben2011changing, ben2011rewarding}. Xianbao et al.~\cite{hu2015behavior} provided information about the travel time for each offered departure time, and alternative time, so the driver does not depend on their experience in choosing the alternatives. 

Recently, \cite{azevedo2018tripod} offered token form incentives for different travel choices such as route, travel modes, and ride-sharing. The proposed model learns individuals' decisions and adapts to their preferences based on their travel history. While these policies were successful in short-term experiments, they do not necessarily result in permanent changes in the travelers' behavior~\cite{kumar2016impacts}.
\cite{papadopoulos2018coordinated} and \cite{kordonis2019mechanisms} provide route and an incentive or fee to volunteer truck drivers to improve the overall traffic condition in a budget balanced mechanism. \cite{papadopoulos2021personalized} considers VOT (Value of Time) in the mechanism to make it personalized.

\vspace{0.2cm}

There have been different choices used as the incentive in transportation studies. \cite{fuji2003does} and \cite{bamberg2003choice} used free bus tickets as an incentive to study the changes in the frequency of traveling by bus. In Australia, an early bird ticket program was offered to relieve the problem of rail overcrowding during peak hours~\cite{currie2011free}. Free WiFi and discounted fares have been effective for Beijing commuters to motivate them to avoid the rush hours of the morning~\cite{zhang2014does}. In the Tripod project~\cite{azevedo2018tripod}, they offered token form incentives that can be redeemed for goods and services at local businesses. In the CAPRI project~\cite{yue2015reducing}, users of the app collected points, and they could trade 100 points for \$1 or use the points to play a game. Knockaert et al.~\cite{knockaert2012spitsmijden} provided smartphones to users who could receive money for their credits or keep the smartphone by collecting enough credits.

\vspace{0.2cm}
 
The preferences of the drivers can be considered in the incentive offering platform. Mohan et al.~\cite{mohan2019exploring} divided effective factors on drivers' decision into two categories of static factors and dynamic factors. The static factors, which are fixed for each person, include the number of available transportation options and the distance of nearby transit. On the other hand, the dynamic factors include weather and travel purpose. They identified the factors from interviews and surveys, and concluded that personalizing can be advantageous for travel assistance. 
Also, the drivers' preferences can be learned through interaction with the individual~\cite{zhu2019personalized, azevedo2018tripod}. Personalized incentives can be used to offer a unique alternative for each individual. The goal is to make the offer close to the individual's preferences and maximize the contribution to the network~\cite{azevedo2018tripod} or to minimize the cost of incentives by maximizing the probability of acceptance of the alternative route~\cite{zhu2019personalized}.
 
\vspace{0.2cm}
 
In this paper, we study the problem of offering personalized incentives to maximize the global utility of the network. Although, previous studies such as \cite{bliemer2009rewarding} consider static rewards for static options like teleworking, biking, and walking, our model assigns different rewards for different alternative routes for different drivers based on the traffic condition and personalization factors. Consequently, we have a much larger system (and a resulting optimization problem).
The implementation of the proposed model could be through a smartphone app where the traffic data can be used to offer incentives to drivers. In addition, smartphones will help the central planner to distribute the computational load for finding the optimal incentive offering strategy. The rest of this paper is organized as follows. In Section~\ref{sec:IncentiveOfferingMechanism}, we present our models. Then to solve our optimization problem efficiently, we propose an efficient distributed algorithm. Results of our numerical experiments are presented in Section~\ref{sec:NumericalExperiments} using data from the Los Angeles area.
\vspace{0.2cm}

\section{Incentive Offering Mechanisms} \label{sec:IncentiveOfferingMechanism}
\vspace{0.2cm}
Let us model the structure of the traffic network with a  directed graph $\cG = (\cV,\cE)$.  
Here $\cV$ is the collection of all major intersections and ramps, which form the set of nodes in the graph. We use the set of edges~$\cE$ to capture the connectivity of the nodes in the graph. Two different nodes are adjacent in the graph if it is possible to directly go from one to another without passing over any other node. The direction of an edge between two nodes is based on the direction of the road from which we can go from one point to another. We also use the notation $|\cE|$ to denote the total number of road segments/edges in our network (i.e. the cardinality of the set~$\cE$). A route is a collection of adjacent edges that starts from one node and ends in another. We use the one-hot encoding scheme to denote the routes. In other words, a given route is represented by a vector $\mathbf{r} \in \{0,1\}^{|\cE|}$. Here, the $k$-th entry of vector~$\mathbf{r}$ is one if the $k$-th edge is a part of route $\mathbf{r}$ and it is zero, otherwise.

\vspace{0.2cm}

Let  $\bT = \{1,\ldots,T\}$ denote the time horizon of interest assuming the system is currently at time $t = 1$. For any $t\in \bT$, we use the random vector $\bv_t \in \mathbb{R}^{|\cE|}$ to represent the traffic volume on different road segments at time~$t$. The $k$-th entry of $\bv_t$ shows the total number of vehicles of road segment~$k$ at time~$t$. Notice that the offered incentives can change the drivers' behavior who are using the platform in the future and thus affecting the vector $\bv_t$. 

\vspace{0.2cm}

We use $\cN$ to denote the set of drivers that we can influence their behavior through offering incentives. For any driver $n \in \cN$, let $\cR_n \subseteq \{0,1\}^{|\cE|}$ denote the set of possible route options for going from its origin to its destination. Let $\cI_n$ be the set of possible incentives we can offer to driver~$n \in \cN$. 
We also use the binary variable $s_i^{\br,n}\in\{0,1\}$ to represent the offered incentives. For any driver~$n\in \cN$ and incentive~$i \in \cI_n$, the variable $s_i^{\br,n} = 1$  if incentive~$i$ is offered to driver~$n$ to take route~$r$; and $s_i^{\br,n} = 0$ otherwise. We assume that we incentivize each driver with only one offer, i.e., $\sum_{\br \in \cR_n} \sum_{i \in \cI_n}s^{\br,n}_i = 1$. Given any incentive offered to the drivers, we model the decision of the drivers stochastically. In particular, we assume after offering incentives, each driver~$n$ chooses route~$\br$ with a certain probability which depends on the amount of incentive, the route, and the driver's preferences, as described below.

\vspace{0.2cm}

The route preferences of drivers depend on different factors such as route travel time, gender, age, and particularly the (monetary) incentive provided to the drivers in our context. Such dependence can be learned using standard machine learning approaches in the presence of data~\cite{xiong2019integrated}. In this project, we rely on the model developed in~\cite{xiong2019integrated} for our preference modeling. We simplify their model by ignoring the less predictive features and only consider two major features: value of incentive and travel time. In particular, we assume that,  given incentive~$i\in \cI_n$ to driver~$n$, the driver chooses route~$\br$ with probability 
\begin{equation} \label{eq:DriverPreference}
p^{\br,n}_{i} = P(\widehat{T}_{\br},i),
\end{equation}
where $\widehat{T}_{\br}$ is the estimate of the travel time for route~$\br$ provided by the incentive offering platform. Notice that when drivers make their routing decisions, they do not know the exact travel time $T_\br$ for route~$\br$, but instead they rely on the estimate~$\widehat{T}_{\br}$ provided by the system. Here, we make an implicit assumption that the drivers do not consider their own judgement about the travel time in their decision. However, if such individual biases for drivers exist, the system can learn it over time  using standard preference learning techniques. Modeling the drivers' behavior in a probabilistic fashion has its own benefits. The decision of a driver for a given incentive amount depends on many factors such as age, gender, income as also studied in~\cite{xiong2019integrated}. It is even likely that the driver's decision may depend on the driver's ``state of mind'' at the time that the incentive is offered. Thus, the features that influence the driver's decision are not completely known to the central planner. In such a setting, probabilistic models can be a better fit for modeling the system. For this reason, in the general area of ``recommendation systems'' in machine learning and statistics, probabilistic models have been widely used to model the behavior of individual users (drivers in our setting)~\cite{zheng2016neural}. 
In addition, we do not assume any traffic control by modeling the probability of drivers' acceptance. Traffic control might be more effective but it needs an authority with the power of changing traffic which is not required in our framework. In our model, drivers can disregard the offers at any time but offers change the probability of accepting drivers' routing choices.

\vspace{0.2cm}

In the next subsections, we present our model and formulation in more detail. For the convenience of the reader, the list of notations defined here and later in the manuscript is presented in Appendix~\ref{apdx:notations}. We present our framework under two different scenarios: First, for simplicity of presenting the ideas, we study the case where it is possible to bring traffic flow below the network capacity. Then, we study the high demand scenario where there is no feasible strategy to bring the demand below the network capacity.

\subsection{Scenario I: Operating Below Network Capacity}
\label{sec:ProposedFormulationI}

Let us first for simplicity assume that there exists a solution that all road segments operate below the capacity. Hence, for that solution, we can assume that the travel time will be based on the free flow traffic. As this section shows, this assumption will result in a mixed integer linear programming optimization which can be solved efficiently using standard solvers.

\vspace{0.2cm}

Given~\eqref{eq:DriverPreference}, the expected value of the volume vector $\mathbf{v}_t$ can be computed as:
\small
\begin{equation} \label{eq:Averagev_t}
    \mathbb{E} \left[\mathbf{v}_t \right]= \sum_{n\in \mathcal{N}} \sum_{i \in \mathcal{I}_n} \sum_{\mathbf{r} \in \mathcal{R}_n} s_i^{\br,n} p^{\mathbf{r},n}_i \boldsymbol{\beta}_{\mathbf{r},t}
\end{equation}
\normalsize
where the vector $\boldsymbol{\beta}_{\br, t}, \in \mathbb{R}^{|\mathcal{E}|}$ shows the probability of being at different links of the network at time $t\in \mathbf{T}$, conditioned on the fact that driver $n$ is on route $\br$. For more details about the vector $\boldsymbol{\beta}_{\br, t}$, please refer to~\cite{ma2018estimating}.

\vspace{0.2cm}

In order to minimize the drivers' total travel time while keeping the volume below the road segment capacity vector~$\mathbf{v}_0$, we need to solve the following optimization problem
\begin{equation} \label{eq:LP2}
\small
\begin{split}
\min_{\{s_i^{\br,n}\}} \quad & \sum_{t\in \bT}\sum_{n\in \mathcal{N}} \sum_{r\in \cR_n} \sum_{i \in \mathcal{I}_n} s_i^{\br,n} p^{\mathbf{r},n}_i \boldsymbol{\beta}_{\mathbf{r},t}^{\intercal} \boldsymbol{\omega} \\
\st \quad & \sum_{n\in \mathcal{N}} \sum_{i \in \mathcal{I}_n} \sum_{\mathbf{r} \in \mathcal{R}_n} s_i^{\br,n} p^{\mathbf{r},n}_i \boldsymbol{\beta}_{\mathbf{r},t}   \leq \bv_0, \quad\forall t  \in \mathbf{T} \\
& \sum_{n\in \mathcal{N}} \sum_{i \in \mathcal{I}_n} \sum_{\br \in \cR_n} s_i^{\br,n} \eta_i \leq \Omega \\
& \sum_{\br \in \cR_n} \sum_{i \in \mathcal{I}_n}s_i^{\br,n} = 1\quad\forall n\in \mathcal{N} \\
& s_i^{\br,n} \in \{0,1\} \quad \forall n\in\mathcal{N}, \forall i\in \mathcal{I}_n, \forall \br \in \cR_n
\end{split}
\normalsize
\end{equation}
where $\boldsymbol{\omega} \in \mathbb{R}^{|\mathcal{E}|}$ is the vector of free flow travel time of the links, $\Omega$ is the total available budget, and $\eta_i$ is the cost of offering incentive~$i$. To keep this optimization problem tractable, we rely on the assumption of a large number of vehicles in each road segment and approximate the random quantity $\bv_t$ with its average $\mathbb{E}[\bv_t]$ provided in equation~\eqref{eq:Averagev_t}. Notice that the objective function is equal to
\begin{equation}
\small
    \min_{\{s_i^{\br,n}\}} \sum_{n\in \mathcal{N}} \sum_{r\in \cR_n} \sum_{i \in \mathcal{I}_n} s_i^{\br,n} p^{\mathbf{r},n}_i \sum_{t\in \bT} \boldsymbol{\beta}_{\mathbf{r},t}^{\intercal} \boldsymbol{\omega} \nonumber
\normalsize
\end{equation}
in which $\sum_{t\in \bT} \boldsymbol{\beta}_{\mathbf{r},t}^{\intercal} \boldsymbol{\omega}$ is the expected travel time of driver $n$ driving on route $\br$.

\vspace{0.2cm}

Problem~\eqref{eq:LP2} is a mixed integer linear program which can be solved via standard  solvers such as Gurobi, AMPL, GAMS, and CPLEX. We use Gurobi in our experiments because of its powerful LP solver.

\subsection{Scenario II:  Operating Above Network Capacity} \label{sec:scenarioII}

In this subsection, we assume that the demand is elevated; thus, there is no incentive offering strategy that can bring the traffic flow below the network capacity. In such a scenario, we still can ``improve'' the congestion via incentivizing individual drivers. Our goal is to optimize a utility of the system as a criterion to compare the traffic condition after incentivizing. To make the formulation more specific, we  use total travel time as the utility function. It is worth noting that while we use this utility, following our steps, one can use any other utility functions such as carbon emissions or energy consumption. 

To compute the total travel time of the system, we sum the travel time of the drivers of all the links over all time periods:
\begin{equation} \label{ObjFunct1}
\small
    \begin{split}
        F_{tt}(\hat{\bv}) = & \sum_{\ell=1}^{|\cE|} \sum_{t=1}^{|\mathbf{T}|} \hat{v}_{\ell, t} \delta_{\ell, t}(\hat{v}_{\ell, t})  \\
    \end{split}
\normalsize
\end{equation}
where $\delta_{\ell, t}$ is the travel time of link $\ell$ at time $t$ (which itself is a function of the volume). Here, $\hat{\bv}$ is the vector of volume of links in which $\hat{v}_{\ell, t}$ is the $(|\cE|\times t + \ell)^{th}$ element of vector $\hat{\bv}$ representing the volume of the $\ell^{th}$ link at time $t$.

\vspace{0.2cm}

To understand the impact of our offered incentives, we estimate the drivers' decision based on the provided incentives, which in turn results in estimating the volume of the links in the horizon of interest. Given these estimated volume values, we  estimate the travel time in the links as described below. 

\vspace{0.2cm}

\noindent\textit{Travel time value $\delta$:} There are different functions that capture the relation between travel time and volume. For example, the link congestion function developed by the Bureau of Public Roads (BPR)~\cite{united1964traffic} defines a nonlinear relation between the volume and travel time of the road segments:   
\begin{equation}
\small
    \begin{aligned}
        f_{\text{BPR}}(v) = t_0 \left(1+0.15\left(\frac{v}{w}\right)^4\right) \nonumber
    \end{aligned}
\normalsize
\end{equation}
where $f_{\text{BPR}}(v)$ is the travel time of the drivers on the link given the assigned traffic volume $v$; the parameter~$t_0$ is the free flow travel time of the link; $v$ is the assigned traffic volume of the link; and $w$ is the practical capacity of the link. We learn $t_0$ and $w$ using historical traffic data. In order to estimate the total travel time of the system, we need to estimate the volume vector $\hat{\bv}$, which we discuss next.

\vspace{0.2cm}

\noindent\textit{Volume vector $\hat{\bv}$:} To compute the volume vector, we need to know the routing  decision of the drivers to be able to (approximately) estimate their location at different times. Clearly, the drivers' decision is a function of the offered incentives. In other words, the location of a driver is dependent on the incentive that we assign to them, because the likelihood of various decisions changes with different incentives. 
Let us first explain our notations for the offered incentives: For each driver, we have a one hot encoded vector describing which route has been incentivized and how much reward has been assigned to it. Thus, for each driver we have a binary vector $\bs_n \in \{0, 1\}^{|\cR|\cdot |\cI|}$ in which only one element has a value of one and it corresponds to the route and the incentive amount that we offer. As we need one vector for each driver, we can aggregate all our incentivization strategy in a matrix $\bS \in \{0, 1\}^{(|\cR|\cdot|\cI|) \times |\cN|}$. Naturally, routes that are not relevant to that OD pair of a driver will get a value of zero in the corresponding incentive vector (since we cannot offer those routes to the driver). 
\vspace{0.2cm}

To understand the drivers' responses to our offered incentives, we need to estimate the probability of acceptance of incentivized routes under different incentives including zero incentive (i.e, no incentive). To model this probability, we use the utility function developed in~\cite{xiong2019integrated} and compute the probability of acceptance of each offered incentive (by using a Softmax function on top of the utility). While the model in~\cite{xiong2019integrated} takes many parameters (such as gender, age, and education of the driver) as input, in our model and numerical experiments we only consider static parameters of the travel time and the reward value  to generate the probability of acceptance of a given incentive/reward.
However, our framework is modular and we can use any prediction model that can estimate drivers' behavior given an incentive amount. We can use any personalized routing model that can learn drivers’ behavior such as a neural network. 
Let $\bP\in [0, 1]^{|\cR| \times (|\cR|\cdot|\cI|)}$ be a matrix encoding the information of probability of picking different routes given the offered (route, incentive) pairs. Thus, the vector $\bP \bS \mathbf{1} \in \mathbb{R}^{|\cR|\times 1}$ shows the expected number of vehicles in each route. 

\vspace{0.2cm}

Given the number of vehicles in each route, the location of each driver for the next time horizon can be modeled in a probabilistic fashion. For this purpose, we rely on the model developed in~\cite{ma2018estimating} where a specific matrix $\bR  \in [0, 1]^{(|\cE|\cdot|\bT|) \times |\cR|} $ is proposed to estimate the probability of the presence of a driver in a given road segment at a specific time in the future (assuming that the driver is picking a specific route). We can compute matrix $\bR$ by running a simulation model if we have enough computation power. In our experiments in subsection~\ref{sec:simModel}, we rely on the historical data in computing matrix $\bR$. Thus, the vector
\begin{equation} \label{xx}
\small
    \begin{split}
        \hat{\bv} =  & \bR \bP \bS \mathbf{1} \in \bbR^{(|\cE|\cdot|\bT|)\times1} \nonumber\\
    \end{split}
\normalsize
\end{equation}
represents the expected number of vehicles in all the links at each time slot. Substituting the expression of $\hat{\bv}$ in \eqref{ObjFunct1}, we get
\begin{equation} \label{eq:F_tt}
\small
    \begin{split}
        F_{tt}(\hat{\bv}) =  & \sum_{\ell=1}^{|\cE|} \sum_{t=1}^{|\mathbf{T}|} (\bA\bS\mathbf{1})_{\ell, t} \delta((\bA\bS\mathbf{1})_{\ell, t}) \\
        = & \sum_{\ell=1}^{|\cE|} \sum_{t=1}^{|\mathbf{T}|} (\ba_{\ell, t}\bS\mathbf{1}) \delta(\ba_{\ell, t}\bS\mathbf{1}) \\
    \end{split}
\normalsize
\end{equation}
where $\ba_{\ell, t}$ is the row of matrix $\bA = \bR \bP$ which corresponds to link $\ell$  at time $t$.
Thus in order to minimize the total travel time in the system via providing incentives to drivers, we need to solve the following optimization problem:
\begin{equation}   
\label{eq:ObjFunctFinal}
\small
    \begin{split}
        \min_{\bS} \quad & \sum_{\ell=1}^{|\cE|} \sum_{t=1}^{|\mathbf{T}|} (\ba_{\ell, t}\bS\mathbf{1}) \delta(\ba_{\ell, t}\bS\mathbf{1}) \\
        \st \quad & \bS^{\intercal} \mathbf{1} = \mathbf{1},\quad 
         \bc^{\intercal} \bS \mathbf{1} \leq \Omega\\
        & \bD \bS \mathbf{1} = \bq,\quad 
         \bS \in \{0, 1\}^{(|\cR||\cI|) \times |\cN|}
    \end{split}
\normalsize
\end{equation}
where $\bc \in \bbR_{+}^{|\cR|\cdot|\cI|}$ is the vector of cost of incentives assigned to each route, $\bD \in \{0, 1\}^{K \times (|\cR|\cdot|\cI|)}$ is the matrix of incentive assignment to the OD pairs, and $\bq \in {\bbR}^{K \times |\cI|}$ is the vector of the number of drivers for each OD pair. Here, $K$ is the number of OD pairs. We explain the constraints in more details below:

\vspace{0.2cm}

\noindent\textbf{Constraint 1} ($\bS^{\intercal} \mathbf{1} = \mathbf{1}$): This constraint simply states that we only  assign  one incentive to each driver. 

\vspace{0.2cm}

\noindent\textbf{Constraint 2} ($\bc^{\intercal} \bS \mathbf{1} \leq \Omega$): This is our budget constraint. The vector $\bc \in \bbR^{|\cR|\cdot|\cI|}$ represents the cost of the different rewards assigned to each driver. $\Omega$ is the total budget.

\vspace{0.2cm}

\noindent\textbf{Constraint 3} ($\bD \bS \mathbf{1} = \bq$): This constraint makes sure that we offer the correct number of rewards for the routes between OD pairs. Recall that $\bS \mathbf{1}$ represents the (expected) number of drivers that have been offered different routes given different rewards. We use matrix $\bD$ to sum the number of drivers that received different reward offers for routes between the same OD pair. $\bq$ is the vector of the actual number of drivers that are travelling between OD pairs and $\bD \bS \mathbf{1}$ must be equal to $\bq$.

\vspace{0.2cm}

\noindent\textbf{Constraint 4} ($\bS \in \{0, 1\}^{(|\cR||\cI|) \times |\cN|}$):  This  constraint imposes binary  structure on our decision parameters. In other words, 0 is not choosing an incentive and 1 is selecting the incentive.

\vspace{0.2cm}

To illustrate our model and the above constraints, we provide an example in Appendix~\ref{appdx:ModelNotationExample}.

\subsection{Algorithm for Offering Incentives and A Distributed Implementation}
\label{sec:ADMMAlgo}
The optimization problem~\eqref{eq:ObjFunctFinal} is of large size while it needs to be solved in almost real time (or  hourly if the drivers send their travel information to the central planner every hour before their trip) in the network. However, due to the existence of binary variable~$\bS$, solving this problem efficiently is difficult\footnote{We conjecture that  problem~\eqref{eq:ObjFunctFinal} is NP-hard to solve  since it is a special instance of polynomial optimization with discrete variables and there does not appear to be any special structure in function $f$ to reduce its complexity.}.
In order to develop an efficient ``approximate'' solver for~\eqref{eq:ObjFunctFinal}, we first relax the binary constraint in~\eqref{eq:ObjFunctFinal} and replace it with the relaxed convex constraint~$ \bS \in [0, 1]^{(|\cR||\cI|) \times |\cN|}$, leading to the relaxed formulation
\begin{equation}
\small
\label{eq:ObjFunctRelaxed}
    \begin{split}
        \min_{\bS} \quad & \sum_{\ell=1}^{|\cE|} \sum_{t=1}^{|\mathbf{T}|} (\ba_{\ell, t}\bS\mathbf{1}) \delta(\ba_{\ell, t}\bS\mathbf{1}) \\
        \st \quad & \bS^{\intercal} \mathbf{1} = \mathbf{1},\quad 
         \bc^{\intercal} \bS \mathbf{1} \leq \Omega\\
        & \bD \bS \mathbf{1} = \bq,\quad
        \bS \in [0, 1]^{(|\cR||\cI|) \times |\cN|}.
    \end{split}
\end{equation}
\normalsize
The constraints in the above optimization problem are clearly convex. By substituting $\ba_{\ell, t}\bS\mathbf{1}$ by $\bga_{\ell, t}$, the objective function becomes a summation of monomial functions with positive coefficients. Moreover, $\bga_{\ell, t}$ is an affine mapping of the optimization variable $\bS$. Since our domain is the nonnegative orthant and monomials are convex in this domain, the objective function is convex. This convexity will allow us to explore the use of standard solvers such as CVX~\cite{cvx}. However, these solvers rely on methods such as interior point methods~\cite{gb08} which requires $O(n^3)$ number of iterations with $n$ being the number of variables. 
This heavy computational complexity prevents us from applying standard solvers. 
In our context, each driver is equipped with a smartphone and; thus, we can distribute the computational burden of solving \eqref{eq:ObjFunctRelaxed} among the drivers. In what follows,  we propose a simple reformulation of the problem leading to a distributed algorithm for solving~\eqref{eq:ObjFunctRelaxed}. To present our algorithm, let us start by reformulating~\eqref{eq:ObjFunctRelaxed} as

\begin{equation}   
\label{eq:ObjFunctRelaxed2}
\small
    \begin{split}
        \min_{\bga,\bu,\bS,\bW,\bH,\bz,\beta} \quad & \sum_{\ell=1}^{|\cE|} \sum_{t=1}^{|\mathbf{T}|} \bga_{\ell, t} \delta(\bga_{\ell, t}) \\ &- \frac{\tilde{\lambda}}{2} \sum_{r=1}^{|\cR|} \sum_{i=1}^{|\cI|}  \sum_{n=1}^{|\cN|}  \bH_{r, i, n}(\bH_{r, i, n}-1)\\
        \st \quad & \bS \mathbf{1} = \bu,\quad
         \bW^{\intercal} \mathbf{1} = \mathbf{1} \\
        & \bD \bu = \bq,  \quad
         \bA \bu = \bga  \quad\quad\quad \\
        & \bH = \bS,  \quad
         \bW = \bS \quad\quad\quad \\
        & \bc^{\intercal} \bu + \beta = \Omega,\quad
          \beta \geq 0  \\
        &  \bH \in [0, 1]^{(|\cR|\cdot|\cI|) \times |\cN|}.\quad\quad\quad\\
    \end{split}
    \normalsize
\end{equation}
\normalsize
This formulation is amenable to the ADMM method~\cite{boyd2011distributed, gabay1976dual, glowinski1975approximation, hong2016convergence, huang2021alternating, barazandeh2021efficient}, which has a natural distributed implementation.
Our ADMM formulation~\eqref{eq:ObjFunctRelaxed2} shows that this computation burden can be distributed among drivers' cell phones. This distributed optimization/federated learning framework can have other standard advantages of federated learning/distributed systems~\cite{li2020federated, lowy2022private}. For example, when proper privacy preserving mechanisms (such as differential privacy~\cite{dwork2006differential}) are utilized, we can guarantee the privacy of drivers since they can participate in the optimization procedure without completely sharing their data.
The steps of this algorithm is summarized in Algorithm~\ref{alg:ADMM-ForOurProblem} and the details of the derivation of its steps is provided in Appendix~\ref{appdx:ADMM}.
\begin{algorithm}[h]
    \caption{Distributed Incentivization via ADMM}
	\begin{algorithmic}[1]
	 \State \textbf{Input}: Initial values: $\bga^0$, $\bS^0$, $\bH^0$, $\bW^0$, $\bu^0$, $\beta^{0}$, $\bla^0_1 \in \mathbb{R}^{|\cR|\cdot|\cI| \times 1}$,  $\bla^0_2 \in \mathbb{R}^{|\cN| \times 1}$,  $\bla^0_3 \in \mathbb{R}^{K \times 1}$,  $\bla^0_4 \in \mathbb{R}^{|\cE| \cdot \bT \times 1}$,  $\bLa^0_5 \in \mathbb{R}^{|\cR|\cdot|\cI| \times |\cN|}$,  $\lambda^0_6 \in \mathbb{R}$,  $\bLa^0_7 \in \mathbb{R}^{|\cR|\cdot|\cI| \times |\cN|}$, Dual update step: $\rho$, Number of iterations: $\tilde{T}$.
        \For {$t = 0, 1, \ldots, \tilde{T}$}
	    	    	    
	    \State $\bu^{t+1} = (\rho\bI +\rho \bD^{\intercal}\bD + \rho\bA^{\intercal}\bA + \rho\bc\bc^{\intercal})^{-1}(\bla_{1}^{t} + \rho\bS^{t}\mathbf{1} - \bD^{\intercal}\bla_3^{t} + \rho\bD^{\intercal}\bq - \bA^{\intercal}\bla_4^{t} + \rho\bA^{\intercal}\bga^{t} - \bc(\lambda_6^{t} + \beta^{t} - \Omega))$  \label{step:u}
	    
	    \State $\bW^{t+1} = (\rho\mathbf{1}\mathbf{1}^{\intercal} + \rho\bI)^{-1}(\rho\mathbf{1}\mathbf{1}^{\intercal} + \rho\bS^{t} - \bLa_7^{t} - \mathbf{1}\bla_2^{t \intercal}) $ \label{step:W}

	   \State $\bH^{t+1} = \mathbbm{1}({\rho>\tilde{\lambda}}) \Pi {\left(\left(\frac{1}{\rho-\tilde{\lambda}}\right)(\rho\bS^{t} - \bLa_5^{t} - \frac{\tilde{\lambda}}{2})\right)}_{[0, 1]} + \mathbbm{1}({\rho<\tilde{\lambda}}) \Pi {\left(\left(\frac{1}{\rho-\tilde{\lambda}}\right)(\rho\bS^{t} - \bLa_5^{t} - \frac{\tilde{\lambda}}{2})\right)}_{\{0, 1\}}$ \label{step:H}

	    \State $\bS^{t+1} = (\rho\bu^{t+1}\mathbf{1}^{\intercal} + \bLa_5^{t} + \rho\bH^{t+1} + \bLa_7^{t} + \rho\bW^{t+1} - \bla_1^{t}\mathbf{1}^{\intercal})(\rho\mathbf{1}\mathbf{1}^{\intercal} + 2\rho\bI)^{-1}$ \label{step:S}
	    
	    \For {$\ell = 0, 1, \ldots, |\cE|$}
	    \For {$\hat{t} = 1, \ldots, |\bT|$}
	    \State $\gamma^{t+1}_{\ell, \hat{t}} = \argmin\limits_{\gamma_{\ell, \hat{t}}}   \gamma_{\ell, \hat{t}} \delta(\gamma_{\ell, \hat{t}})  + \bla^{t}_{4, (\ell, \hat{t})}(\ba_{\ell, \hat{t}} \bu^{t} - \gamma_{\ell, \hat{t}}) + \frac{\rho}{2} (\ba_{\ell, \hat{t}} \bu^{t} - \gamma_{\ell, \hat{t}})^{2}$ \label{step:gamma}
	    \EndFor
	    \EndFor
	    
	    \State $\beta^{t+1} = \Pi\left(\Omega - \bc^{\intercal}\bu^{t+1} -\frac{1}{\rho}\lambda_6^{t} \right)_{\bbR_{+}}$ \label{step:beta}

	    \State $\bla^{t+1}_{1} = \bla^{t}_{1} + \rho (\bS^{t+1}\mathbf{1} - \bu^{t+1})$  \label{step:lambda1}
	    	    
	    \State $\bla^{t+1}_{2} = \bla^{t}_{2} + \rho (\bW^{t+1 \intercal} \mathbf{1} -\mathbf{1})$ \label{step:lambda2}
	    	    
	    \State $\bla^{t+1}_{3} = \bla^{t}_{3} + \rho (\bD \bu^{t+1} - \bq)$ \label{step:lambda3}
	    	    
	    \State $\bla^{t+1}_{4} = \bla^{t}_{4} + \rho (\bA \bu^{t+1} - \bga^{t+1})$ \label{step:lambda4}
	    	    
	    \State $\bLa^{t+1}_{5} = \bLa^{t}_{5} + \rho (\bH^{t+1}  - \bS^{t+1})$ \label{step:lambda5}
	    	    
	    \State $\lambda^{t+1}_{6} = \lambda^{t}_{6} + \rho (\bc^{\intercal} \bu^{t+1} + \beta^{t+1} - \Omega)$ \label{step:lambda6}
	    	    
	    \State $\bLa^{t+1}_{7} = \bLa^{t}_{7} + \rho (\bW^{t+1} - \bS^{t+1})$ \label{step:lambda7}
	    
		\EndFor
		\State \textbf{Return:} $\bS^{\tilde{T}}$
	\end{algorithmic}
\label{alg:ADMM-ForOurProblem}
\end{algorithm}

\vspace{0.2cm}

In Algorithm~\ref{alg:ADMM-ForOurProblem}, $\Pi(\cdot)_{[0,1]}$ is the projection operator that projects each entry of the input matrix to the interval $[0,1]$ and $\Pi(\cdot)_{\bbR_{+}}$ is the projection operator that projects each entry of the input matrix on to $\bbR_{+}$. For a better convergence of Algorithm~\ref{alg:ADMM-ForOurProblem}, we randomly permute the order of update of two blocks of ADMM at each iteration. The first block consists of variables $\bu, \bW$, and $\bH$ and the second block includes $\bS, \gamma_{\ell, \hat{t}}$, and $\beta$. In Algorithm~\ref{alg:ADMM-ForOurProblem} the steps~\ref{step:W}, \ref{step:H}, and \ref{step:S} are computationally cumbersome due to the size of the matrices $\bW,\bH,$ and $\bS$. However, notice that each column of the matrices $\bW,\bH$, and $\bS$ corresponds to single driver and hence the computation corresponding to each column can be performed in parallel on smartphone devices of the drivers. Moreover, since the steps are not coupled, they can be solved in parallel on the driver's smart devices. Theorem~\ref{thm:ConvergenceADMM} guarantees the convergence of our ADMM algorithm.

\begin{theorem} \label{thm:ConvergenceADMM}
Algorithm~\ref{alg:ADMM-ForOurProblem} finds an $\epsilon$-optimal solution of problem~\eqref{eq:ObjFunctRelaxed} in $O(1/\epsilon)$ iterations~\cite{he20121}.
\end{theorem}

\vspace{0.2cm}

In optimization problem~\eqref{eq:ObjFunctRelaxed} (and consequently~\eqref{eq:ObjFunctRelaxed2}), all solutions $\bS^{\ast}$ with a fixed value of $\bS^{\ast} \mathbf{1}=\bu^{\ast}$ lead to the same objective as long as $\bS^{\ast \intercal} \mathbf{1}= \mathbf{1}$. Hence, this convex problem can have infinite number of solutions (with many of them not even close to binary). Therefore, in order to find (approximately) binary solutions, we add the following regularizer to the objective function in~\eqref{eq:ObjFunctRelaxed2}:
\small
\begin{equation}
    \Re(\bH_{r, i , n}) = - \frac{\tilde{\lambda}}{2} \bH_{r, i , n}(\bH_{r, i , n}-1)
\end{equation}
\normalsize
where $\tilde{\lambda} \in \bbR_{+}$ is the regularization parameter and $\bH_{r, i, n} \in [0, 1]$. This regularizer forces the elements of matrix $\bH$ to be as close as possible to the binary domain $\{0, 1\}$.

While Algorithm~\ref{alg:ADMM-ForOurProblem} returns the solution of the optimization problem~\eqref{eq:ObjFunctRelaxed2}, this problem~\eqref{eq:ObjFunctRelaxed2} is a relaxation of the original problem~\eqref{eq:ObjFunctFinal}. Hence, the obtained solution in Algorithm~\ref{alg:ADMM-ForOurProblem} must be utilized to obtain a feasible point in~$\eqref{eq:ObjFunctFinal}$. For this step, we solve the following mixed integer (linear) problem
\begin{equation}   
\label{eq:linearModelS}
\small
    \begin{split}
        \min_{\bS} \quad & \| \bS\mathbf{1} - \bu^{\ast}\|_{1} \\
        \st \quad & \bS^{\intercal} \mathbf{1} = \mathbf{1},\quad
         \bc^{\intercal} \bS \mathbf{1} \leq \Omega\\
        & \bD \bS \mathbf{1} = \bq, \quad 
         \bS \in \{0, 1\}^{(|\cR||\cI|) \times |\cN|}
    \end{split}
\normalsize
\end{equation}
where $\bu^{\ast}$ is the optimal solution obtained by Algorithm~\ref{alg:ADMM-ForOurProblem}. We can use off-the-shelf solvers such as Gurobi to solve~\eqref{eq:linearModelS}.

The BPR function of Algorithm~\ref{alg:ADMM-ForOurProblem} in Section~\ref{sec:ProposedFormulationI} can capture both Scenario~I in Section~\ref{sec:ProposedFormulationI} and Scenario~II in Section~\ref{sec:scenarioII}. However, the computational requirements for the free-flow case in Scenario~I in Section~\ref{sec:ProposedFormulationI} is less expensive compared to that of the congested case in Scenario~II in Section~\ref{sec:scenarioII}. Thus, the model~\eqref{eq:LP2} in Scenario~I in Section~\ref{sec:ProposedFormulationI} is an alternative when computational resources are limited.

\vspace{0.2cm}
\section{Numerical Experiments} \label{sec:NumericalExperiments}
We evaluate the performance of our algorithms using data from the Los Angeles area. The Los Angeles region is ideally suited for being the validation area because there are multiple routes connecting most OD pairs. Additionally, researchers at the University of Southern California have developed the Archived Data Management System (ADMS) that collects, archives, and integrates a variety of transportation datasets from Los Angeles, Orange, San Bernardino, Riverside, and Ventura Counties. ADMS includes access to real-time traffic data from 9500 highway and arterial loop detectors with measurements every 30 seconds and 1 minute respectively. 

\vspace{0.2cm}

Due to the lack of access to the drivers' routing information, we need to estimate the origin-destination (OD) matrix from the network flow information. Rows and columns of the OD matrix correspond to the origin and destination points respectively. For OD matrix $A$, the element $A_{(i,j)}$ is the number of drivers going from point $i$ to point $j$. The OD matrix estimation problem is severely under-determined~\cite{van1980most, cascetta1984estimation, bell1991real}. There are two categories of OD matrices: static and dynamic~\cite{bera2011estimation}. Due to the high resolution of our data, most of the existing dynamic OD estimation (DODE) methods become computationally inefficient. In addition, we do not have prior data of the OD matrix which many studies consider as given data~\cite{krishnakumari2019data, carrese2017dynamic, nigro2018exploiting, kim2014using} and we do not have access to prior observations of the OD matrix. Given these barriers, we relied on the algorithm proposed by~\cite{ma2018estimating}. This algorithm performs without employing any prior OD matrix information.

\subsection{Simulation Model}
\label{sec:simModel}
In our numerical experiments, we integrate different datasets and models to evaluate the performance of our algorithms. First, we extract the speed data, volume data, and sensor information including the location of sensors from the Archived Data Management System (ADMS). Then, we use the distances of sensors, extracted from the location of sensors using Google Maps API, to create the graph of the network. We created three sets of graph networks corresponding to the regions depicted in Fig.~\ref{fig:region_x4}, Fig.~\ref{fig:region_y2}, and Fig.~\ref{fig:region_y3}. In the next step, the speed data, volume data, and the network graph are used for estimating OD pairs by the algorithm provided in~\cite{ma2018estimating}. The total number of estimated incoming drivers for all three experiments is presented in Fig.~\ref{fig:number_driver}. For each OD pair, we find up to 4 different routing options. In particular, we start by the shortest path for each OD pair. Then, we remove the edges in this path and go with the second shortest path, and we continue this process until we find 4 different routes between the origin and destination (or no other routes exist). We use the model in~\cite{xiong2019integrated} to compute the acceptance probability for the different offers on the different routes for each individual driver. The parameters used for computing the probability are static values provided by~\cite{xiong2019integrated} 
and we only calibrate some of the parameters because we are only using historical data and we do not have access to the drivers' features such as age or gender.
In this paper, we do not learn the route choice model of the drivers so the parameters of the probability model are fixed but it is possible to adapt our routing model to the drivers’ preferences by observing the drivers’ behavior. 
We run three different experiments that model the road network at different scales.  In Experiment~I, we model an arterial region (Fig.~\ref{fig:region_x4}) but includes surface streets. For Experiment~II, we model a large  network of highways (Fig.~\ref{fig:region_y2}). For Experiment~III, we model a moderate region (Fig.~\ref{fig:region_y3}) which is a subset of the region in Experiment~II. 
We use travel time savings as our metric for performance evaluation in all experiments. Beside travel time savings, we also include the monetary value of traffic reduction based on the Value of Time (VOT) as an alternative metric. Our base Value of Time (VOT) is derived from the estimation of \cite{castillo2020benefits} which is \$2.63 per minute or \$157.8 per hour.
As we discussed in section~\ref{sec:ProposedFormulationI}, model~\eqref{eq:LP2} assumes that there exists a traffic flow solution operating below the network capacity. When this assumption is not satisfied, our model results in an ``approximate'' solution. To evaluate the validity of this approximation in heavily congested networks, we run model~\eqref{eq:LP2} for heavily congested networks (Experiments~I and II) with many alternative routes so that we can reasonably reduce the congestion level. As we see in Experiments~I and II, this model can provide a reasonable approximation in both arterial (Experiment~I) and highways (Experiment~II) and leads to congestion reduction even when the final result is above the system capacity. Experiment~III, on the other hand, shows a limitation of model~\eqref{eq:LP2} where incentives offered does not necessarily translate to reduced congestion. This is because of the heavy congestion and the lack of availability of enough alternative routes to reduce congestion (so that the final solution is far away from the free flow traffic and a linear approximation of travel time is no longer accurate enough). In this scenario, our numerical experiments demonstrate the superiority of model~\eqref{eq:ObjFunctFinal} in reducing congestion even in such a heavily congested network. We were not able to run model~\eqref{eq:ObjFunctFinal} for Experiments~I and II due to the large number of nodes in the network. However, relying on edge computation, this model could be solved efficiently in practice as we discussed in subsection~\ref{sec:ADMMAlgo}. 
To solve model~\eqref{eq:LP2}, we use the Gurobi solver in all experiments. Also, we solve model~\eqref{eq:ObjFunctFinal} in Experiment~III utilizing Gurobi and MOSEK to compare their results with Algorithm~\ref{alg:ADMM-ForOurProblem}. 
Gurobi and MOSEK are the state of the art off-the-shelf commercial solvers of linear and mixed integer optimization problems.
To have a better balance between accuracy and required time for solving the problem, we set the relative mixed integer programming optimality gap at 0.01 for both Gurobi and MOSEK in the experiments.
The comparison between ADMM, Gurobi, and MOSEK is shown in Experiment~III.
While we only provide incentives to the drivers that enter the system in the first time interval,  our incentive offering mechanism considers estimations of the traffic flows in the next time intervals. 
The selected drivers for incentivization are from the same cohort. We randomly select a group of drivers between 7 AM and 7:15 AM. Then, we use the selected drivers to compare the performance of the model with different budget values on the total travel time for 7 AM to 8 AM. 
While our formulation is static, it 
can be applied in a dynamic environment if 
solved frequently in the network in order to offer incentives to the drivers. 

\subsection{Experiment~I}
\begin{figure}[] 
  \centering
  \begin{tabular}{@{}c@{}}
    \includegraphics[width=0.4\linewidth]{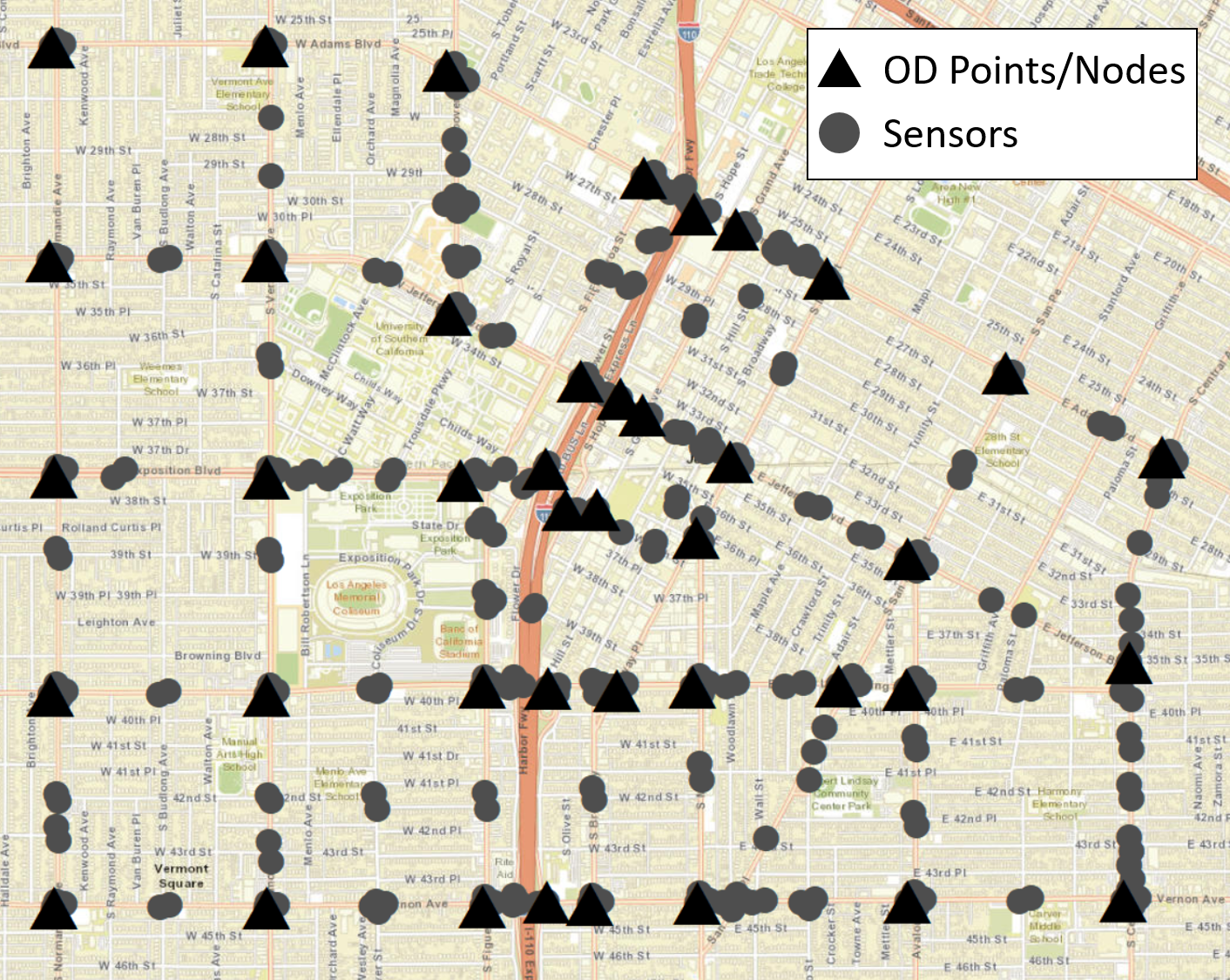} \\[\abovecaptionskip]
  \end{tabular}
  \vspace{-0.25cm}
  \caption{Studied region in Experiment~I.}
  \label{fig:region_x4}
\end{figure}

In Experiment~I, we check the performance of model~\eqref{eq:LP2} using the ADMS data for May 5\textsuperscript{th}, 2018 with the incentive set $\cI=\{\$0,\$2,\$10\}$. The studied region, which is depicted in Fig.~\ref{fig:region_x4}, includes the data of 301 sensors. Based on the ADMS data, we created a graph with 41 nodes, 139 links, and $105.5$ miles of road. OD points are located on intersections and close to the ramp of the highways. The number of OD pairs is 1681 and there are 4278 paths between them in total. We assume 7494 drivers enter the system between 7 AM and 8 AM and we consider 1805 drivers entering the system in the first 15 minutes for incentivization. 
Results of model~\eqref{eq:LP2} at 100\% penetration rate are presented in TABLE~\ref{table:lp2_region_x4_78AM}. In this table, ``Total travel time'' shows the travel time computed via the BPR function after offering incentives. When the budget is increased from \$1000 to \$10,000, the percentage of total travel time decrease is improved by almost a factor of 2 from $4.07\%$ to $7.37\%$. The row ``Cost'' in the table shows the amount of the budget that was used. In all cases, almost all the budget is used. The value of saved time is much larger than budget except in the budget of \$10000. Note that a budget of zero is the case of no incentive.
\begin{table}[]
\small
\footnotesize
\centering
\begin{tabular}{c|c|c|c|c|}
\cline{2-5}
                        & \multicolumn{4}{c|}{Budget (\$)}                         \\ \cline{2-5} 
                       & $0$     & $100$         & $1000$         & $10000$\\ \hline
\multicolumn{1}{|c|}{\makecell{Cost (\$)}} & $0$     & $100$         & $1000$         & $9996$\\ \hline
\multicolumn{1}{|c|}{\makecell{Value of saved\\ time (\$)}} & $0$     & $775$         & $4930$         & $8927$\\ \hline
\multicolumn{1}{|c|}{\makecell{Total travel\\ time (hour)}} & $767$    & $762$         & $736$         & $710$\\ \hline
\end{tabular}
\vspace{+0.15cm}
\caption{Experiment~I: Linear model~\eqref{eq:LP2}.}
\label{table:lp2_region_x4_78AM}
\vspace{-0.3cm}
\normalsize
\end{table}

Notice that model~\eqref{eq:LP2} may result in an infeasible optimization problem (particularly in a heavily congested network). Hence, we included a parameter $\alpha$ in our formulation as the multiplier of the allowed capacity. In other words, we use $\alpha\times \bv_0$ instead of $\bv_0$ in model~\eqref{eq:LP2}. We only consider this multiplier during the computation of incentives; however, during the computation of total travel time, we use the original capacity. We assumed zero dollar incentive in our probabilistic model for drivers that are not receiving an incentive.
\begin{table}[]
\small
\footnotesize
\centering
\begin{tabular}{|c|c|c|c|c|}
\hline
                  \makecell{Number\\ of drivers\\ entering\\ the system} & \makecell{Budget (\$)} & \makecell{\% of \\ rewarded\\drivers} & \makecell{Average\\ incentive \\amount} 
                  & \makecell{Reduction\\ in total\\ travel time} \\ \hline
\multirow{2}{*}{7402} & 1000 & 6.67\% & \$2.00 & 4.07\% \\ \cline{2-5} 
                  & 10000 & 13.64\% & \$9.78 & 7.37\%\\ \hline
\end{tabular}
\vspace{+0.25cm}
\caption{Comparison of \$1000 and \$10000 budget in Experiment~I.}
\label{table:comprasion_results}
\vspace{-0.2cm}
\normalsize
\end{table}

TABLE~\ref{table:comprasion_results} shows that increasing the budget results in higher percentage of drivers to whom we offered the incentive and a higher average amount of offered incentives. In addition, we observe that even offering incentives to 6.67\% of the drivers (with an average of \$2.00  monetary incentive per driver) can reduce the total travel time by 4.07\%. If approximately 13.64\% of the drivers are incentivized with an average of \$9.78 per driver, the total travel time can be reduced by almost 7.37\%. For more details about the distribution of offered incentives to drivers in Experiment~I, please see TABLE~\ref{table:distIncentiveI_7-8} in the Appendix. 

Fig.~\ref{fig:region_x4_LP_penet_TT} shows the effect of the penetration rate (percentage of drivers who are considered in incentivization) on travel time decrease. By reducing the penetration rate, we experience smaller travel time decrease because of the flexibility of the model in selecting drivers decreases. Although reducing the penetration rate adversely affects the incentivization, the model focuses on available drivers for reducing the travel time. For more details on the numbers provided in Fig.~\ref{fig:region_x4_LP_penet_TT}, please see TABLE~\ref{table:PenEffectEx1} and TABLE~\ref{table:PenEffectHourEx1} in the Appendix.

\begin{figure}[] 
  \centering
  \begin{tabular}{@{}c@{}}
    \includegraphics[width=0.5\linewidth]{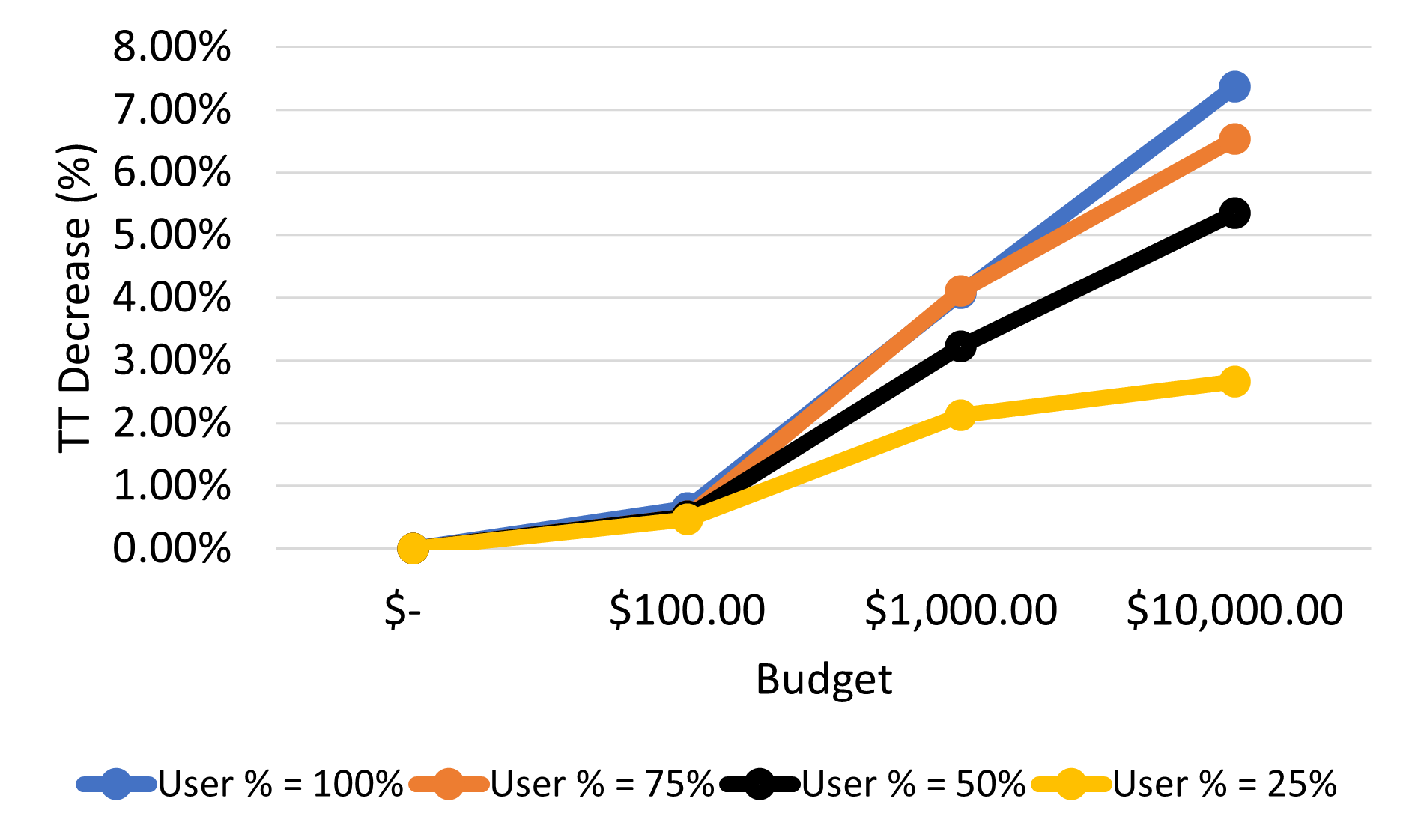} \\[\abovecaptionskip]
  \end{tabular}
  \vspace{-0.25cm}
  \caption{Effect of the penetration rate on percentage of travel time decrease in Experiment~I.}
  \label{fig:region_x4_LP_penet_TT}
\end{figure}

\subsection{Experiment~II}
\begin{figure}[] 
  \centering
  \begin{tabular}{@{}c@{}}
    \includegraphics[width=0.4\linewidth]{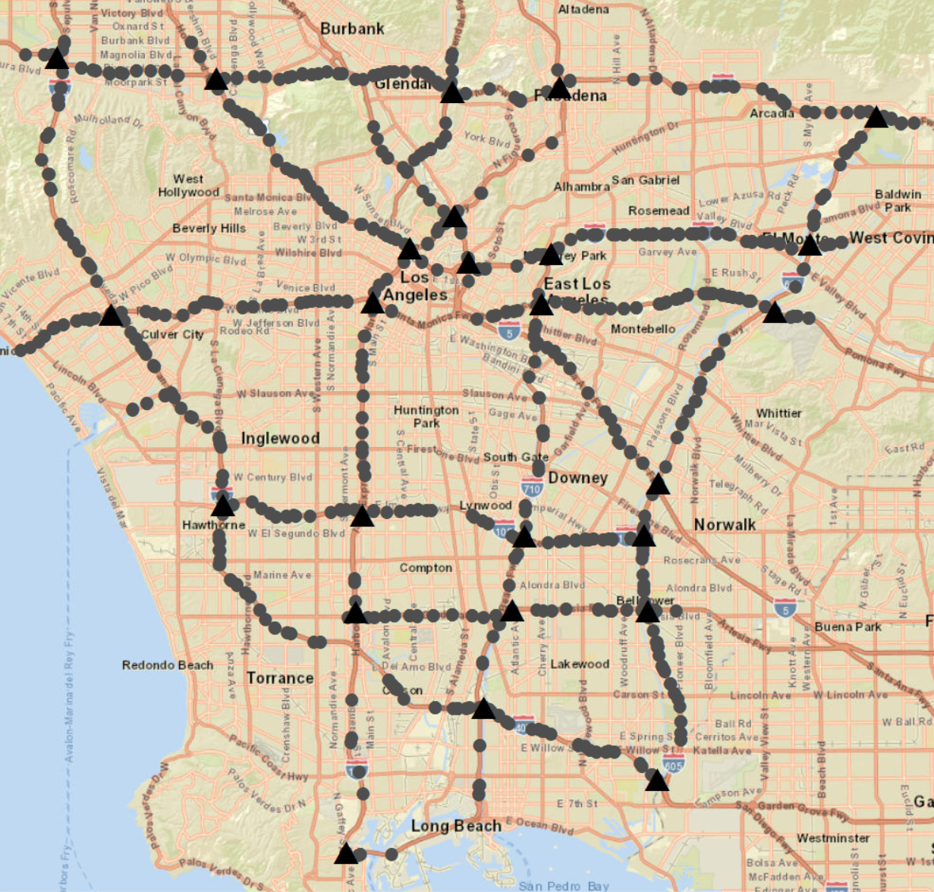} \\[\abovecaptionskip]
  \end{tabular}  
  \vspace{-0.35cm}
  \caption{Studied region in Experiment~II.}
  \label{fig:region_y2}
\end{figure}

In Experiment~II, we evaluate the performance of our methods for the region depicted in Fig.~\ref{fig:region_y2} with 753 sensors under two different possible sets of incentives:
\begin{enumerate}
    \item[$\bullet$] $\cI_{1}=\{\$0, \$2,\$10\}$
    \item[$\bullet$] $\cI_2=\{\$0,\$1,\$2,\$3,\$5,\$10\}$
\end{enumerate}
This region only includes data of highway sensors with 25 OD points and 32 links which includes $707.6$ miles of road. The number of OD pairs is 625, and there are 1331 paths between them in total. We assume 15093 drivers enter the system between 7 AM and 8 AM. Our incentivization model considers 4126 drivers entering the system in the first 15 minutes. 
The results of our experiment at 100\% penetration rate are presented in TABLE~\ref{table:lp2_region_y2_set1_78AM} for incentive set $\cI_1$, and in TABLE~\ref{table:lp2_region_y2_set2_78AM} for incentive set $\cI_2$. The value of saved time is much larger than budget for all budget values for both incentive set and it can go up to 97 times of the cost.
\begin{table}[]
\small
\footnotesize
\centering
\begin{tabular}{c|c|c|c|c|}
\cline{2-5}
                        & \multicolumn{4}{c|}{Budget (\$)}                         \\ \cline{2-5} 
                      & $0$     & $100$         & $1000$         & $10000$\\ \hline
\multicolumn{1}{|c|}{\makecell{Cost (\$)}} & $0$     & $100$         & $1000$         & $9998$\\ \hline
\multicolumn{1}{|c|}{\makecell{Value of saved\\ time (\$)}} & $0$     & $6378$         & $12919$         & $15712$\\ \hline
\multicolumn{1}{|c|}{\makecell{Total travel\\ time (hour)}} & $5023$    & $4983$   & $4941$       & $4924$\\ \hline
\end{tabular}
\vspace{+0.25cm}
\caption{Experiment~II: Linear model~\eqref{eq:LP2} for incentive set $\cI_1$.}
\label{table:lp2_region_y2_set1_78AM}
\vspace{-0.3cm}
\normalsize
\end{table}

\begin{table}[]
\small
\footnotesize
\centering
\begin{tabular}{c|c|c|c|c|}
\cline{2-5}
                        & \multicolumn{4}{c|}{Budget (\$)}                         \\ \cline{2-5} 
                      & $0$     & $100$         & $1000$         & $10000$\\ \hline
\multicolumn{1}{|c|}{\makecell{Cost (\$)}} & $0$     & $100$         & $1000$         & $10000$\\ \hline
\multicolumn{1}{|c|}{\makecell{Value of saved\\ time (\$)}} & $0$     & $9791$         & $18499$        & $28535$\\ \hline
\multicolumn{1}{|c|}{\makecell{Total travel\\ time (hour)}} & $5023$    & $4961$         & $4906$         & $4843$\\ \hline
\end{tabular}
\vspace{+0.25cm}
\caption{Experiment~II: Linear model~\eqref{eq:LP2} for incentive set $\cI_2$.}
\label{table:lp2_region_y2_set2_78AM}
\normalsize
\vspace{-0.5cm}
\end{table}

\begin{table}[]
\small
\centering
\footnotesize
\begin{tabular}{c|c|c|c|c|c|}
\cline{2-6}
                                        &        \makecell{Number\\ of drivers\\ entering\\ the system}           & \makecell{Budget\\ (\$)}  & \makecell{\% of \\ rewarded\\drivers} & \makecell{Average\\ incentive \\amount} & \makecell{Reduction\\ in total\\ travel time} \\ \hline
\multicolumn{1}{|c|}{\multirow{2}{*}{\makecell{Exp. II-1\\}}} & \multirow{2}{*}{\makecell{15093\\}} & 1000 & 2.97\% & \$2.23 & 1.63\% \\ \cline{3-6} 
\multicolumn{1}{|c|}{}   &    & 10000 & 16.42\% & \$4.03 & 1.98\% \\ \hline

\multicolumn{1}{|c|}{\multirow{2}{*}{\makecell{Exp. II-2}}} & \multirow{2}{*}{15093} & 1000 & 3.87\% & \$1.71 & 2.33\% \\ \cline{3-6} 
\multicolumn{1}{|c|}{}   &   & 10000 & 15.97\% & \$4.15 & 3.60\% \\ \hline

\end{tabular}
\caption{Comparison of \$1000 and \$10000 budget in Experiment~II.}
\label{table:comprasion_results2}
\normalsize
\end{table}

In addition to confirming our previous observations in Experiment~I, Experiment~II shows the  diversity gain related to the incentive set $\cI_1$  (see the ``Reduction in total travel time'' column in TABLE~\ref{table:comprasion_results2}). In other words, more choices in the incentive set provides more flexibility for the algorithm, resulting in the total travel time reduction. 
For more details about the distribution of offered incentives to drivers in Experiment~II, please see TABLE~\ref{table:distIncentiveII_settingI_7-8} and TABLE~\ref{table:distIncentiveII_settingII_7-8} in the Appendix. By examining Experiments~I and II, we observe that more alternative routes leads to more gain in travel time reduction.

\subsection{Experiment~III}
\begin{figure}[] 
  \centering
  \begin{tabular}{@{}c@{}}
    \includegraphics[width=0.5\linewidth]{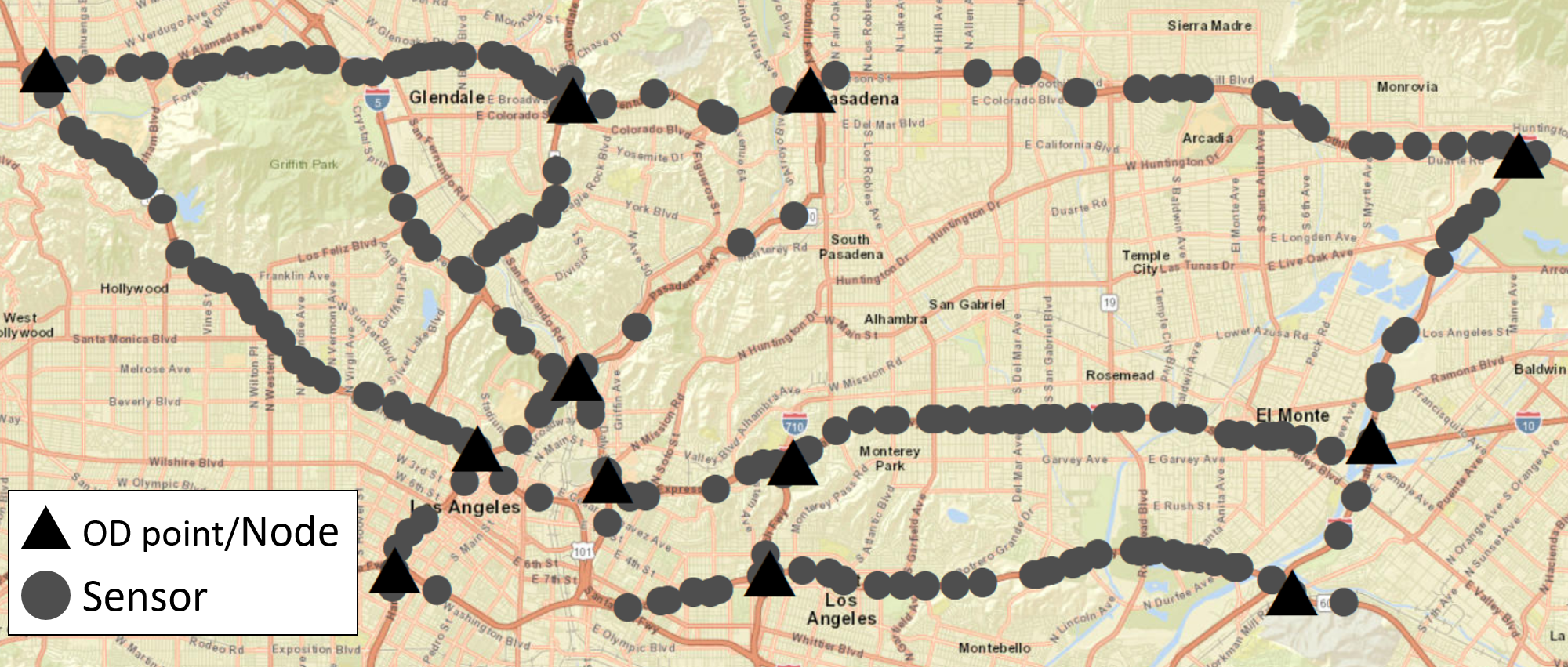} \\[\abovecaptionskip]
  \end{tabular}
  \vspace{-0.35cm}
  \caption{Studied region in Experiment~III.}
  \label{fig:region_y3}
\end{figure}

\begin{figure}[] 
  \centering
  \begin{tabular}{@{}c@{}}
    \includegraphics[width=0.5\linewidth]{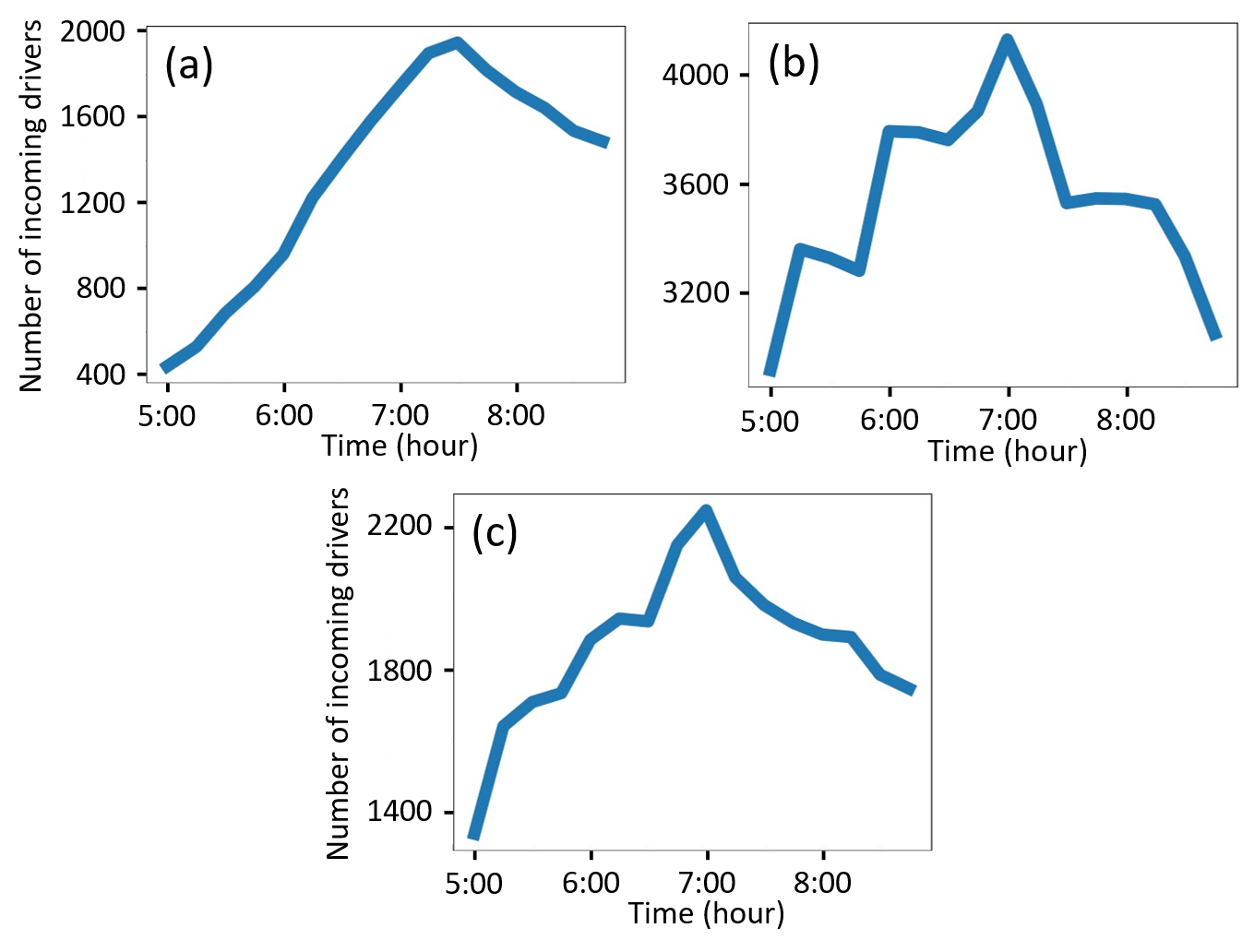} \\[\abovecaptionskip]
  \end{tabular}
  \vspace{-0.7cm}
  \caption{Total estimated number of drivers entering the system (in 15-minute intervals). (a) Experiment~I, (b) Experiment~II, and (c) Experiment~III.}
  \label{fig:number_driver}
\end{figure}

In Experiment~III, we compare the performance of linear model~\eqref{eq:LP2} and ADMM model~\eqref{eq:ObjFunctFinal} using the incentive set $\cI=\{\$0,\$2,\$10\}$. The region considered in our analysis is depicted in Fig.~\ref{fig:region_y3}. This region includes the data of 293 sensors. Based on the ADMS data, we created a graph with 12 nodes, 32 links, and $288.1$ miles of road. The number of OD pairs is 144 and there are 270 paths between them in total. The estimated total number of drivers incoming to the system between 5 AM to 9 AM by the OD estimation algorithm is depicted in Fig.~\ref{fig:number_driver} (c). 
In our simulations, we assume 8220 drivers enter the system between 7 AM and 8 AM. Our incentivization model considers 2248 drivers entering the system in the first 15 minutes.
Results of model~\eqref{eq:LP2} and model~\eqref{eq:ObjFunctFinal}  at 100\% penetration rate are presented in TABLE~\ref{table:resutls_region_y3_78AM}. 

\begin{table}[]
\small
\centering
\footnotesize
\begin{tabular}{c m{4em}|c|c|c|c|}
\cline{3-6}
                                        &  & \multicolumn{4}{c|}{Budget (\$)} \\ \cline{3-6} 
                                        &  &  $0$   &  $100$   &  $1000$   &  $10000$   \\ \hline
\multicolumn{1}{|c|}{\multirow{3}{*}{\makecell{Model\\\eqref{eq:LP2}\\ (Linear)}}} & Cost (\$) & $0$ &  $100$   &  $1000$   &   $10000$  \\ \cline{2-6} 
\multicolumn{1}{|c|}{}                  &  \makecell{Value of\\saved\\ time (\$)} & $0$     & $-7265$         & $1498$         & $-25118$\\ \cline{2-6} 
\multicolumn{1}{|c|}{}                  & Total travel time (hour) &  3005.22   &  3051   &  2995  &  3164\\ \hline
\multicolumn{1}{|c|}{\multirow{2}{*}{\makecell{Model\\\eqref{eq:ObjFunctFinal}\\ (ADMM)}}} & Cost (\$) &  $0$   & $100$     &  $1000$   &   $9954$  \\ \cline{2-6} 
\multicolumn{1}{|c|}{}                  &  \makecell{Value of\\saved\\ time (\$)} & $0$     & $4105$         & $37016$         & $53646$\\ \cline{2-6} 
\multicolumn{1}{|c|}{}                   & Total travel time (hour) &  3005   &  2979   &  2771   &  2665   \\ \hline
\end{tabular}
\vspace{+0.25cm}
\caption{Experiment~III: Linear model~\eqref{eq:LP2} and model~\eqref{eq:ObjFunctFinal}.}
\label{table:resutls_region_y3_78AM}
\normalsize
\vspace{-0.4cm}
\end{table}

\begin{figure}[!htb]
   \begin{minipage}{0.48\textwidth}
     \centering
     \includegraphics[width=1.0\linewidth]{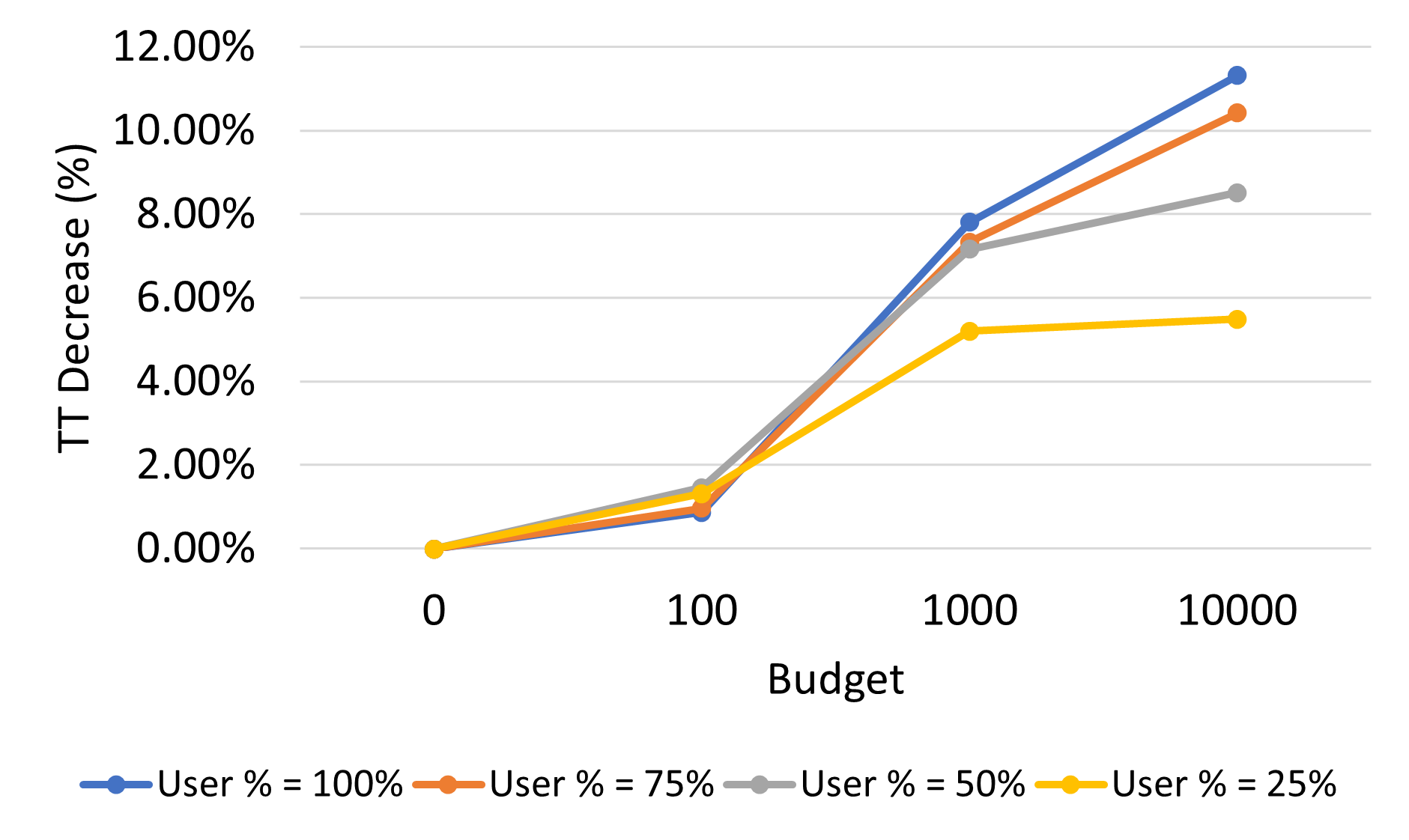}
     \caption{Effect of the penetration rate on percentage of travel time decrease in Experiment~III, model~\eqref{eq:ObjFunctFinal}, Algorithm~\ref{alg:ADMM-ForOurProblem}.}\label{fig:region_y3_ADMM_penet_TT}
   \end{minipage}\hfill
   \begin{minipage}{0.48\textwidth}
     \centering
     \includegraphics[width=1.0\linewidth]{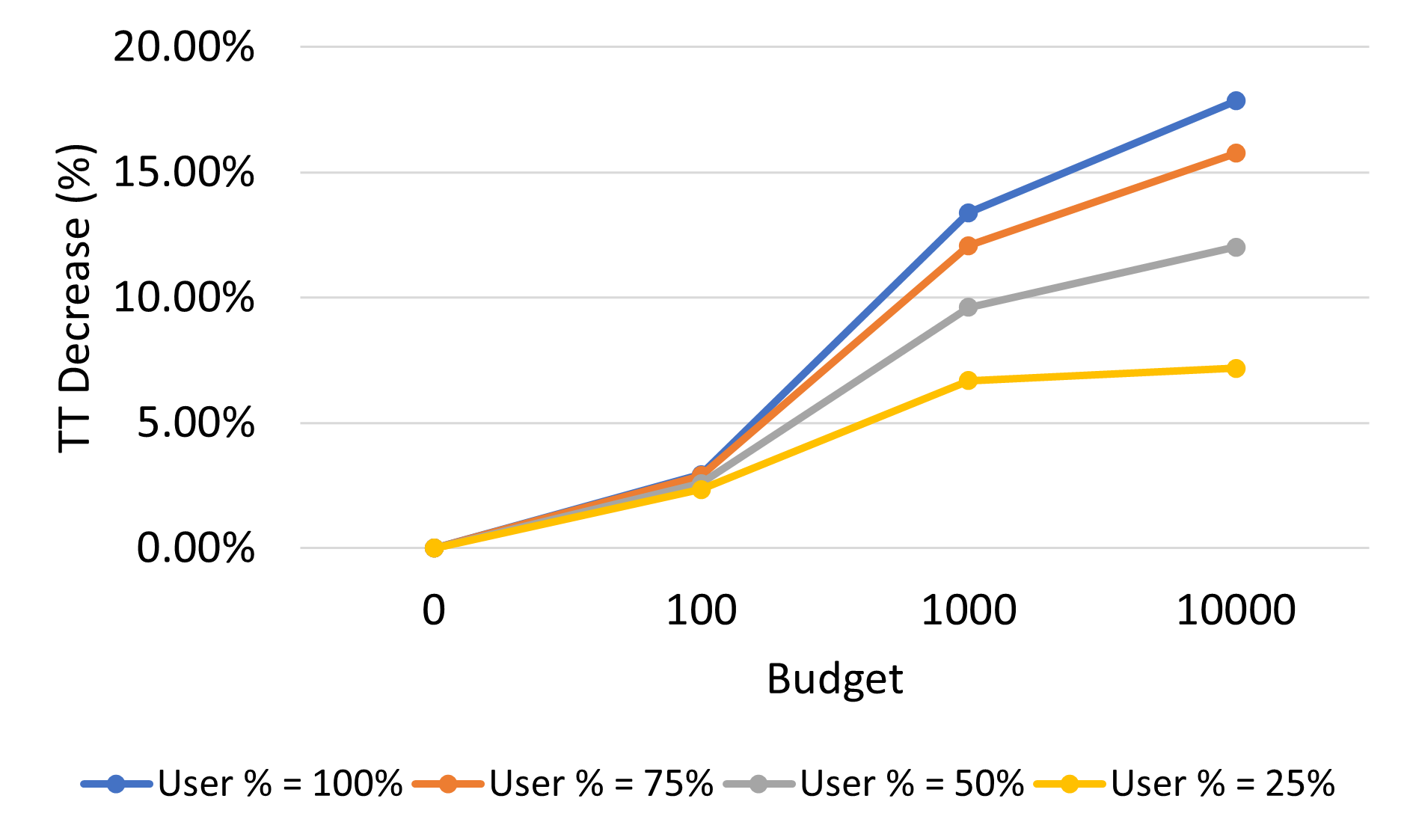}
     \caption{Effect of the penetration rate on percentage of travel time decrease in Experiment~III, model~\eqref{eq:ObjFunctFinal}, Gurobi solver.}\label{fig:region_y3_Gurobi_penet_TT}
   \end{minipage}
\end{figure}

\begin{table}[]
\small
\centering
\footnotesize
\begin{tabular}{c|c|c|c|c|c|}
\cline{2-6}
                                        &        \makecell{Number\\ of drivers\\ entering\\ the system}           & \makecell{Budget\\ (\$)}  & \makecell{\% of \\ rewarded\\drivers} & \makecell{Average\\ incentive \\amount} & \makecell{Reduction\\ in total\\ travel time} \\ \hline
\multicolumn{1}{|c|}{\multirow{2}{*}{\makecell{Exp.~III\\Model~\eqref{eq:LP2}}}} & \multirow{2}{*}{\makecell{8220\\}} & 1000 & 6.08\% & \$2.00 & 0.32\% \\ \cline{3-6} 
\multicolumn{1}{|c|}{}                  &                   & 10000 & 15.86\% & \$7.71 & -5.30\% \\ \hline

\multicolumn{1}{|c|}{\multirow{2}{*}{\makecell{Exp. III\\Alg.~\ref{alg:ADMM-ForOurProblem}}}} & \multirow{2}{*}{8220} & 1000 & 5.26\% & \$2.31 & 7.81\% \\ \cline{3-6} 
\multicolumn{1}{|c|}{}                  &                   & 10000 & 16.17\% & \$7.49 & 11.31\% \\ \hline

\end{tabular}
\caption{Comparison of \$1000 and \$10000 budget in Experiment~III.}
\label{table:comprasion_results3}
\normalsize
\end{table}
TABLE~\ref{table:comprasion_results3} shows a major shortcoming of model~\eqref{eq:LP2}. Although the objective function in model~\eqref{eq:LP2} reduces the total free flow travel time, the actual travel time is not based on the free flow travel time. Hence, model~\eqref{eq:LP2} does not necessarily decrease the total travel time when the available budget increases as can be seen in TABLE~\ref{table:comprasion_results3}. However, the decrease in total travel time is guaranteed in model~\eqref{eq:ObjFunctFinal} for larger budget using Gurobi, MOSEK, and Algorithm~\ref{alg:ADMM-ForOurProblem}. Model~\eqref{eq:ObjFunctFinal} directly minimizes the travel time based on the BPR function so it captures the nonlinear relation between travel time and volume. Moreover, when the volume is greater than the capacity, model~\eqref{eq:ObjFunctFinal} which is a more accurate model representing the traffic network produces better results. 
These phenomena can be observed in Fig.~\ref{fig:region_y3_ADMM_penet_TT} and Fig.~\ref{fig:region_y3_Gurobi_penet_TT}. We have not reported the travel time decrease by MOSEK because its results are less than $0.6\%$ different with the results of Gurobi. Although both Gurobi and MOSEK find slightly better solutions for model~\eqref{eq:ObjFunctFinal}, it can take take up to 35 hours for Gurobi and up to 80 hours for MOSEK to solve the problem. However, Algorithm~\ref{alg:ADMM-ForOurProblem} takes at most 1.12 hours. Algorithm~\ref{alg:ADMM-ForOurProblem} by offering incentives to 16.17\% of the vehicles at 100\% penetration rate in Experiment~III (with an average of \$7.49 monetary incentive per driver)  can reduce the total travel time by 11.31\% using model~\eqref{eq:ObjFunctFinal}. For a budget of $\$1000$, model~\eqref{eq:ObjFunctFinal} has 7.92\% larger reduction in percentage of travel time compared to model~\eqref{eq:LP2} although both offer almost the same amount of incentive on average to the almost same percentage of drivers. For a \$10000 budget, not only model~\eqref{eq:LP2} does not compete with model~\eqref{eq:ObjFunctFinal} but also it increases the total travel time. Despite this shortcoming of model~\eqref{eq:LP2}, it is still useful in low congested areas and when the computational resources are limited. In most cases, model~\eqref{eq:LP2} found a reasonable approximate solution except when the budget was \$10000.  
The computation time of model~\eqref{eq:LP2} is 2.6 minutes, but model~\eqref{eq:ObjFunctFinal} requires up to 14.12 hours to run if we employ  serial computation. Utilizing parallel computation as described in section~\ref{sec:ADMMAlgo}, we can reduce the computational time to at most 1.12 hours. 
If we decompose matrices A and D beforehand among the drivers in our parallel computation, we can further reduce the computational time to 0.11 hours (or 6 minutes). The value of saved time using Algorithm~\ref{alg:ADMM-ForOurProblem} is much larger than budget for all budget values and it can go up to 41 times of the cost. 
For more details about the distribution of the offered incentives to the drivers in Experiment~III, please see TABLE~\ref{table:distIncentiveIII_Linear_7-8}, TABLE~\ref{table:distIncentiveIII_ADMM_7-8}, and TABLE~\ref{table:distIncentiveIII_Gurobi_7-8} in the Appendix. 
The effect of the penetration rate on travel time decrease in Experiment~III for model~\eqref{eq:ObjFunctFinal} is depicted in Fig.~\ref{fig:region_y3_ADMM_penet_TT} and Fig.~\ref{fig:region_y3_Gurobi_penet_TT}. The behavior is similar to our observation in Fig.~\ref{fig:region_x4_LP_penet_TT} which was for Model~\eqref{eq:LP2}. For more details of the numbers provided in Fig.~\ref{fig:region_y3_ADMM_penet_TT} and Fig.~\ref{fig:region_y3_Gurobi_penet_TT}, please see TABLE~\ref{table:PenEffectEx3ADMM}, TABLE~\ref{table:PenEffectHourEx3ADMM}, TABLE~\ref{table:PenEffectEx3Gurobi}, and TABLE~\ref{table:PenEffectHourEx3Gurobi} in the Appendix.

\section{Conclusion}
In this paper, we developed mathematical models and proposed algorithms for offering personalized incentives to drivers to reduce congestion in the network. In this framework, drivers share their origin-destination and routing information with a central planner. Based on this information, the central planner then offers incentives to drivers to incentivize/enforce a socially optimal routing strategy. The incentives are offered based on solving large-scale optimization problems in our framework. In our framework, we bring together prior works to model the behavior of drivers in response to the offered incentives as well as the resulting congestion reduction in the network where no traffic control is required. 
We paid special attention to minimizing the total travel time of the network. In addition, we showed that this problem can be solved in a distributed fashion where some of the computations are performed on individual drivers’ smart devices. Finally, we evaluated the performance of our models and algorithms using Archived Data Management System (ADMS) data. Our experiments showed that the proposed framework can lead up to a 11\% decrease in the total travel time of the system during rush hour times.

In this work, the incentives are only offered to alter the routing decision of the drivers. In future work, it is crucial to look at the effect of offering incentives to change the mode or time of the drivers’ trips. These options will bring additional flexibility to the model, which in turn will result in further congestion reduction. To compute the drivers' acceptance probability, we can include more aspects of drivers' characteristics and features into account such as gender, age. and salary. In addition, we can utilize preference learning in computing drivers' acceptance probability if we have access to the data of preferences of the drivers. All the codes for this project can be found in~\cite{MyCode}.

\bibliographystyle{unsrtnat}
\bibliography{bibliography}

\section{List of Notations} \label{sec:notation}
\label{apdx:notations}
The following symbols are used in this paper. 
\begin{itemize}
    \item $\cG$: Directed graph of the traffic network
    \item $\cV$: Set of nodes of graph $\cG$ which correspond to major intersections and ramps
    \item $\cE$: Set of edges of graph $\cG$ which correspond to the set of road segments
    \item $|\mathcal{E}|$: Total number of road segments/edges in the network $\cG$ (i.e. the cardinality of the set $\mathcal{E}$)
    \item $\mathbf{r}$: Route vector 
    \item $\mathbf{T}$: Time horizon
    \item $|\mathbf{T}|$: Number of time units (i.e. the cardinality of $\mathbf{T}$)
    \item $\mathbf{v}_0$: Capacity vector of road segments
    \item $\mathbf{v}_t$: Volume vector of road segments at time $t$ 
    \item $\mathcal{N}$: Set of drivers
    \item $|\mathcal{N}|$: Number of drivers (i.e. the cardinality of the set $\mathcal{N}$)
    \item $\mathcal{R}_n$: Set of possible route options for driver $n$
    \item $\mathcal{R}$: Total set of possible route options for all drivers
    \item $|\mathcal{R}|$: Number of possible route options (i.e. the cardinality of the set $\mathcal{R}$)
    \item $\mathcal{I}_n$: Set of possible incentives to offer to driver $n$
    \item $\mathcal{I}$: Total set of possible incentives to all drivers
    \item $|\mathcal{I}|$: Number of possible incentives (i.e. the cardinality of the set $\mathcal{I}$)
    \item $s_i^{\mathbf{r},n}$: Decision parameter indicates whether incentive $i$ is offered to driver $n$ for route $\br$
    \item $p_{i}^{\mathbf{r},n}$: The probability of acceptance of route $\mathbf{r}$ by driver $n$ given incentive $i$
    \item $\widehat{T}_{\br}$: The estimate of the travel time for route $\br$ provided by the incentive offering platform
    \item $T_{\br}$: The exact travel time for route $\br$ 
    \item $\boldsymbol{\beta}_{\mathbf{r},t}$: The vector of location of driver that is traveling a route $\br$ at time $t$
    \item $\eta_i$: The cost of incentive $i$
    \item $F_{tt}(.)$: Total travel time function
    \item $\delta_{\ell, t}$: Travel time of link $\ell$ at time $t$
    \item $\hat{\bv}$: The vector of volume of links at different times in the horizon
    \item $\hat{v}_{\ell, t}$: The $(|\cE|\times t + \ell)^{th}$ element of vector $\hat{\bv}$ representing the volume of $\ell^{th}$ link at time $t$
    \item $t_0$: The free flow travel time of the link
    \item $v$: The traffic volume of the link
    \item $w$: The practical capacity of the link
    \item $\bs_n$: The binary decision vector for one driver in which only one element has value of one and it corresponds to the route and the incentive amount that we offer
    \item $f_{BPR}(.)$: BPR function
    \item $\bS$: Decision matrix
    \item $\bR$: The matrix of location of a driver 
    \item $\bP$: Route choice probability matrix
    \item $\bD$: The matrix of incentive assignment to OD pairs
    \item $\bq$: The vector of number of drivers for each OD pair
    \item $\bc$: The vector of cost of incentives assigned to each route
    \item $\Omega$: Budget
    \item $\boldsymbol{\omega}$: The vector of free flow travel time of links
    \item $\ba_{\ell, t}$: The row of matrix $\bA = \bR \bP$ which corresponds to link $\ell$  at time $t$
    \item $K$: The number of OD pairs
    \item $e$: An edge of graph $\cG$ which corresponds to a road segments in the traffic network
\end{itemize}

\vspace{0.2cm}

\section{Details of Alternating Direction Method of Multipliers (ADMM)}
\label{appdx:ADMM}
In this section, we explain the details of solving~\eqref{eq:ObjFunctRelaxed} utilizing ADMM method~\cite{boyd2011distributed}.
Let
\small
\begin{align}
    & \mathcal{L}(\bga, \bS, \bH, \bW, \bu, \beta)  \nonumber\\
    \triangleq & \; F_{tt}(\bga) + \mathbb{I}_{[0, 1]^{(|\cR|\cdot|\cI|) \times |\cN|}}(\bH)+\mathbb{I}_{\mathbb{R}^{+}}(\beta) \nonumber\\
    &  + \langle\bla_1, \bS\mathbf{1}-\bu\rangle + \langle\bla_2, \bW^{\intercal}\mathbf{1} - \mathbf{1}\rangle + \langle\bla_3, \bD\bu-\bq\rangle \nonumber\\
    &  + \langle\bla_4, \bA\bu - \bga\rangle +  \langle\bLa_5, \bH - \bS\rangle \nonumber\\
    &   + \lambda_6(\bc^{\intercal}\bu  + \beta - \Omega) + \langle\bLa_7, \bW - \bS\rangle    \nonumber\\
    & + \frac{\rho}{2}||\bS\mathbf{1} - \bu||^{2} + \frac{\rho}{2}||\bW^{\intercal}\mathbf{1}  -\mathbf{1}||^{2} \nonumber\\
    &  + \frac{\rho}{2}||\bD\bu - \bq||^{2}  + \frac{\rho}{2}||\bA\bu - \bga||^{2}    \nonumber\\
    & + \frac{\rho}{2}||\bH - \bS||^{2} +  \frac{\rho}{2}(\bc^{\intercal}\bu + \beta - \Omega)^{2} \nonumber\\
    &    + \frac{\rho}{2}||\bW - \bS||^{2} - \frac{\tilde{\lambda}}{2} \sum_{\ell=1}^{|\cE|} \sum_{t=1}^{|\mathbf{T}|}  \bH_{\ell, t}(\bH_{\ell, t}-1) \nonumber\\
\end{align}
\normalsize
be the augmented Lagrangian function of \eqref{eq:ObjFunctRelaxed2} with the set of Lagrange multipliers $\{\bla_1, \bla_2, \dots, \bLa_7\}$ and $\rho > 0$ be the primal penalty parameter. Then, ADMM solves \eqref{eq:ObjFunctRelaxed2} by the following iterative scheme

\small
\begingroup
\allowdisplaybreaks
\begin{align*}
\bu^{t+1} = & \argmin\limits_{\bu} \quad  \langle\bla_1^{t}, \bS^{t+1}\mathbf{1} - \bu\rangle +  \langle\bla_3^{t}, \bD\bu - \bq\rangle   \\
&  + \langle\bla_4^{t}, \bA\bu-\bga\rangle + \lambda_6(\bc^{\intercal}\bu  + \beta - \Omega)  \\
& + \frac{\rho}{2}||\bS^{t+1}\mathbf{1} - \bu||^{2} + \frac{\rho}{2}||\bD\bu - \bq||^2 \\
& + \frac{\rho}{2}||\bA\bu - \bga||^2  +  \frac{\rho}{2}(\bc^{\intercal}\bu + \beta - \Omega)^{2}\\
\bW^{t+1} = & \argmin\limits_{\bW} \quad  \langle\bla_2^{t}, \bW^{\intercal}\mathbf{1} - \mathbf{1}\rangle + \langle\bLa_7^{t}, \bW - \bS^{t+1}\rangle \\
& + \frac{\rho}{2}||\bW^{\intercal}\mathbf{1} - \mathbf{1}||^2  +  \frac{\rho}{2}||\bW - \bS^{t+1}||^2  \\
&\\
\bH^{t+1} = & \argmin\limits_{\bH}  \quad \mathbbm{1}({\rho>\tilde{\lambda}})\mathbb{I}_{[0, 1]^{(|\cR|\cdot|\cI|) \times |\cN|}}(\bH)\\
& + \mathbbm{1}({\rho<\tilde{\lambda}})\mathbb{I}_{\{0, 1\}^{(|\cR|\cdot|\cI|) \times |\cN|}}(\bH) \nonumber\\
& + \langle\bLa_5^{t}, \bH- \bS^{t+1}\rangle + \frac{\rho}{2}||\bH-\bS^{t+1}||^{2} \nonumber\\
& - \frac{\tilde{\lambda}}{2} \sum_{\ell=1}^{|\cE|} \sum_{t=1}^{|\mathbf{T}|}  \bH_{\ell, t}(\bH_{\ell, t}-1) \\
\bS^{t+1} = & \argmin\limits_{\bS}  \quad \langle\bla_1^{t}, \bS\mathbf{1}-\bu^{t}\rangle + \langle\bLa_5^{t}, \bH^{t} - \bS\rangle      \\
&  + \langle\bLa_7^{t}, \bW^{t} - \bS\rangle +  \frac{\rho}{2}||\bS\mathbf{1} - \bu^{t}||^{2}  \\
& + \frac{\rho}{2}||\bH^{t}- \bS||^{2} + \frac{\rho}{2}||\bW^{t} - \bS||^{2}  \\
\beta^{t+1} = & \argmin\limits_{\beta} \quad  \mathbb{I}_{+}(\beta) + \lambda_6^{t}\left(\bc^{\intercal}\bu^{t+1} + \beta - \Omega\right) \\
& + \frac{\rho}{2}\left(\bc^{\intercal}\bu^{t+1} + \beta - \Omega\right)^2  \\
 \bla^{t+1}_{1} = & \bla^{t}_{1} + \rho \left(\bS^{t+1}\mathbf{1} - \bu^{t+1}\right)   \\
\bla^{t+1}_{2} = & \bla^{t}_{2} + \rho \left(\bW^{t+1 \intercal} \mathbf{1} -\mathbf{1}\right)  \\
\bla^{t+1}_{3} = & \bla^{t}_{3} + \rho \left(\bD \bu^{t+1} - \bq\right)  \\
\bla^{t+1}_{4} = & \bla^{t}_{4} + \rho \left(\bA \bu^{t+1} - \bga^{t+1}\right)  \\
\bLa^{t+1}_{5} = & \bLa^{t}_{5} + \rho \left(\bH^{t+1}  - \bS^{t+1}\right)  \\
\lambda^{t+1}_{6} = & \lambda^{t}_{6} + \rho \left(\bc^{\intercal} \bu^{t+1} + \beta^{t+1} - \Omega\right)  \\
\bLa^{t+1}_{7} = & \bLa^{t}_{7} + \rho \left(\bW^{t+1} - \bS^{t+1}\right)   
\end{align*}
\endgroup
\normalsize

The primal update rules can be simplified as
\small
\begingroup
\allowdisplaybreaks
\small
\begin{align*}
\gamma^{t+1}_{\ell, \hat{t}} = & \argmin\limits_{\gamma_{\ell, \hat{t}}}   \gamma_{\ell, \hat{t}} \delta(\gamma_{\ell, \hat{t}}) + \bla^{t}_{4, (\ell, \hat{t})}(\ba_{\ell, \hat{t}} \bu^{t} - \gamma_{\ell, \hat{t}})\\
& + \frac{\rho}{2} (\ba_{\ell, \hat{t}} \bu^{t} - \bga_{\ell, \hat{t}})^{2}, \quad \forall \ell, \forall \hat{t}  \\
\bS^{t+1} = &\; \frac{1}{\rho} (-\bla_1^{t} \mathbf{1}^{\intercal} + \bLa_{5}^{t} + \bLa_{7}^{t} + \rho \bu^{t} \mathbf{1}^{\intercal} + \rho \bH^{t}\\
& + \rho \bW^{t})(\mathbf{1}\mathbf{1}^{\intercal} + 2\bI)^{-1}  \\
\bH^{t+1} =  &\; \mathbbm{1}(\rho>\tilde{\lambda})\Pi {\left(\left(\frac{1}{\rho-\tilde{\lambda}}\right)(\rho\bS^{t} - \bLa_5^{t} - \frac{\tilde{\lambda}}{2})\right)}_{[0, 1]}\\ & +\mathbbm{1}(\rho<\tilde{\lambda})\Pi {\left(\left(\frac{1}{\rho-\tilde{\lambda}}\right)(\rho\bS^{t} - \bLa_5^{t} - \frac{\tilde{\lambda}}{2})\right)}_{\{0, 1\}}  \\
\bW^{t+1} = &\; \frac{1}{\rho} (\bI + \mathbf{1}\mathbf{1}^{\intercal})^{-1} (-\mathbf{1}\bla_2^{t \intercal} - \bLa_7^{t} + \rho \mathbf{1}\mathbf{1}^{\intercal} + \rho \bS^{t+1})  \\
\bu^{t+1} = &\; \frac{1}{\rho} (\bI + \bD^{\intercal} \bD + \bA^{\intercal}\bA + \bc \bc^{\intercal})^{-1} (\bla_1^{t} - \bD^{\intercal} \bla_3^{t} - \bA^{\intercal} \bla_4^{t}\\
& + \rho\bS^{t+1}\mathbf{1} + \rho \bD^{\intercal}\bq + \rho \bA^{\intercal}\bga^{t+1} - \lambda_6 \bc - \beta \rho\bc + \Omega\rho \bc)\\
\beta^{t+1} = & \; \Pi\left( \frac{1}{\rho} (-\lambda_6^{t} -\rho \bc^{t \intercal} \bu^{t} + \rho \Omega)\right)_{\bbR_{+}}  
\end{align*}
\normalsize
\endgroup
\normalsize

\vspace{0.2cm}

\section{An Example of the Model and Notations}
\label{appdx:ModelNotationExample}

In this section, we present a small example of a network to illustrate our model and notations. Consider the  network
\begin{figure}[H]
  \centering
  \begin{tabular}{@{}c@{}}
    \includegraphics[width=.3\linewidth]{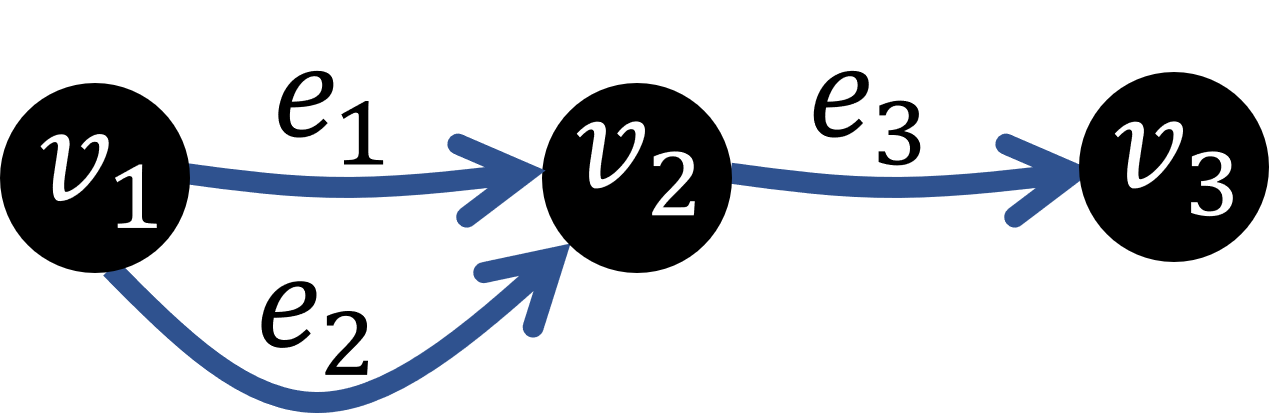} \\[\abovecaptionskip]
  \end{tabular}
  \caption{Network example $\cG_1$.}\label{fig:network}
\end{figure}
\noindent
where $\cV = \{\nu_1, \nu_2, \nu_3\}$ is the set of nodes and $\cE = \{e_1, e_2, e_3\}$ is the set edges (roads). Details of the links and attributes are represented in TABLE~\ref{fig:tableLinks}. The (origin, destination) pair is ($\nu_1$, $\nu_3$). There are two routes going from origin to destination as illustrated in TABLE~\ref{fig:tableRoutes}. The time horizon set is $\bT = \{1, 2, 3\}$ and each time is $0.2$ hour. To estimate the location of drivers at each time, we need matrix $\bR \in [0, 1]^{9\times6}$ as follows
\vspace{-0.3cm}

\small
\begin{figure}[H]
\small
\begin{align*}
\bR = 
\begin{blockarray}{ccccccc}
\small
& \makecell{t_1\\=1\\ \br_1} & \makecell{t_1\\=1\\ \br_2} & \makecell{t_1\\=2\\ \br_1} & \makecell{t_1\\=2\\ \br_2} & \makecell{t_1\\=3\\ \br_1} & \makecell{t_1\\=3\\ \br_2} \\
\begin{block}{c(cccccc)}
  \makecell{t_2=1, e_1} & 1 & 1 & 0 & 0 & 0 & 0 \\
  \makecell{t_2=1, e_2} & 0 & 0 & 0 & 0 & 0 & 0 \\
  \makecell{t_2=1, e_3} & 0.5 & 0 & 0 & 0 & 0 & 0 \\
  \makecell{t_2=2, e_1} & 0 & 0 & 1 & 1 & 0 & 0 \\
  \makecell{t_2=2, e_2} & 0 & 1 & 0 & 0 & 0 & 0 \\
  \makecell{t_2=2, e_3} & 0.5 & 0 & 0.5 & 0 & 0 & 0 \\
  \makecell{t_2=3, e_1} & 0 & 0 & 0 & 0 & 1 & 1 \\
  \makecell{t_2=3, e_2} & 0 & 0 & 0 & 1 & 0 & 0 \\
  \makecell{t_2=3, e_3} & 0 & 0 & 0.5 & 0 & 0.5 & 0 \\
\end{block}
\end{blockarray}
\end{align*}
\end{figure}
\normalsize
\noindent
where $t_1$ is the entrance time of the driver and $t_2$ is the driver's arrival time at the road. In model \eqref{eq:LP2}, the column vector $\boldsymbol{\beta_{\br, t}}$ corresponds to the columns of matrix $\bR.$\\
\small
\begin{table}[]
\small
\centering
\begin{tabular}{c|c|c|c|}
\cline{2-4}
                       &  \makecell{Length\\ (Mile)} & \makecell{Speed\\ (mph)} & \makecell{Travel time\\ (Hour)} \\ \hline
\multicolumn{1}{|c|}{$e_1$} &  5 & 50 & 0.1\\ \hline
\multicolumn{1}{|c|}{$e_2$} &  10 & 50 & 0.2\\ \hline
\multicolumn{1}{|c|}{$e_3$} &  5 & 50 & 0.1\\ \hline
\end{tabular}
\vspace{+0.15cm}
\caption{Set of edges.} \label{fig:tableLinks}
\end{table}
\normalsize

\vspace{-0.4cm}

\small
\begin{table}[]
\small
\centering
\small
\begin{tabular}{c|c|c|}
\cline{2-3}
                       &  $\br$ & Graph\\ \hline
\multicolumn{1}{|c|}{\makecell{Route 1\\ $e_1 \rightarrow e_3$}} &  $\br_1=\begin{bmatrix} 1 \\ 0 \\ 1\end{bmatrix}$ & \includegraphics[width=.2\linewidth]{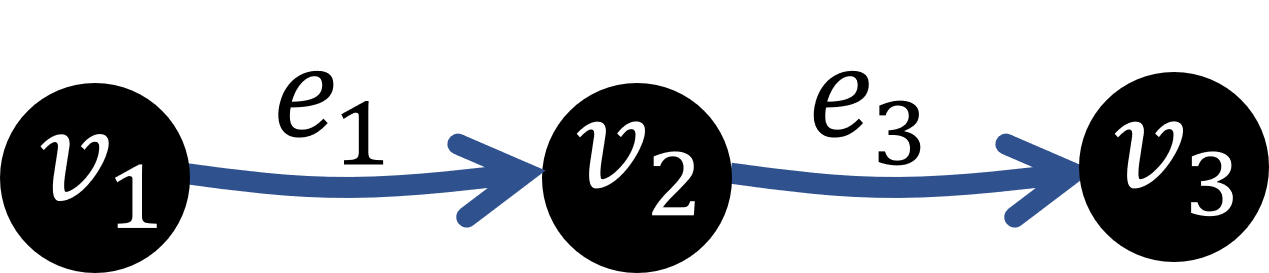} \\ \hline
\multicolumn{1}{|c|}{\makecell{Route 2\\ $e_2 \rightarrow e_3$}} & $\br_2 =\begin{bmatrix} 0 \\ 1 \\ 1\end{bmatrix}$ & \includegraphics[width=.2\linewidth]{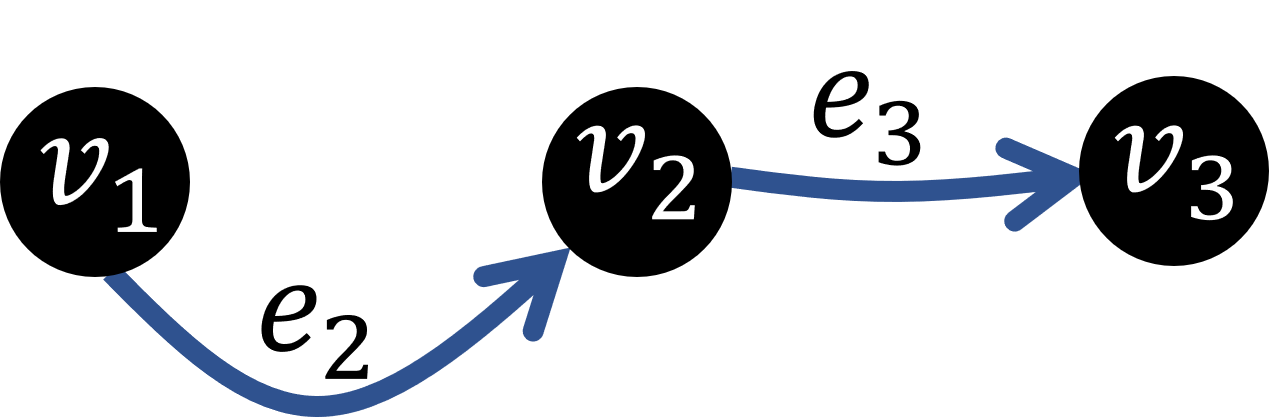} \\ \hline
\end{tabular}
\vspace{+0.15cm}
\caption{Set of routes.} \label{fig:tableRoutes}
\end{table}
\normalsize
\vspace{-0.4cm}

Assume there are two  drivers in the system and $\cN = \{d_1, d_2\}$. We want to offer rewards from the set $\cI = \{\$0, \$5\}$ to control the traffic. To estimate the probability of choosing routes given an offered incentive at a time, we use matrix $\bP \in [0, 1]^{6\times12}$ when incentive $i$ is offered:
\vspace{-0.5cm}
\begin{figure}[H]
\small
\begin{align*}
\bP_{t_i} = 
\begin{blockarray}{ccccc}
& \multicolumn{2}{c}{\;\;\;\text{No incentive}} & \$5 \rightarrow \br_1 & \$5 \rightarrow \br_2  \\
\begin{block}{c(cccc)}
   t=1, \br_1 & \;\;\; 0.50 & 0.50 & 0.97 & 0.03 \\
   t=1, \br_2 & \;\;\; 0.50 & 0.50 & 0.03 & 0.97 \\
   t=2, \br_1 & \;\;\; 0.50 & 0.50 & 0.97 & 0.03 \\
   t=2, \br_2 & \;\;\; 0.50 & 0.50 & 0.03  & 0.97 \\
   t=3, \br_1 & \;\;\; 0.50 & 0.50 & 0.97 & 0.03 \\
   t=3, \br_2 & \;\;\; 0.50 & 0.50 & 0.03  & 0.97 \\
\end{block}
\end{blockarray}
\\,\forall i \in \{1, 2, 3\}
\end{align*}
\end{figure}
\vspace{-0.55cm}
\small
\begin{align*}
\bP = \begin{bmatrix} 
 \bP_{t_1} &  \bP_{t_2} &  \bP_{t_3} \\
\end{bmatrix}
\end{align*}
\normalsize

Probability matrices for all three times are equal because the speed is the same in all three times. We compute the probability of choosing route $k$ given that $\$i'$ is offered for route $j'$ by
\begin{equation}
\small
\begin{split}
    & \mathbb{P}(\br=k, i=(\$i' \rightarrow \text{route} \; j')) \\
     = & \frac{\exp{( -0.086 tt_k + 0.7 i' \mathbb{I}_{k=j'})}}{\exp{(-0.086tt_{j'} + 0.7 i')} + \sum_{j\neq j'} \exp{(-0.086tt_j)}}
\end{split}
\normalsize
\end{equation}
\noindent
where $tt_j$ is the travel time of route $j$. We use~\cite{xiong2019integrated} to extract these coefficients.

\vspace{0.3cm}

\section{Details of the Numerical Experiments}
\label{appdx:detailsNumerical}

\begin{table}[H]
\centering
\begin{tabular}{cc|c|ccc|}
\cline{3-6}
\multicolumn{2}{c|}{\multirow{2}{*}{}}                                                                                                          & \multirow{2}{*}{Budget} & \multicolumn{3}{c|}{Incentive}                              \\ \cline{4-6} 
\multicolumn{2}{c|}{}                                                                                                                           &                         & \multicolumn{1}{c|}{\$0}  & \multicolumn{1}{c|}{\$2} & \$10 \\ \hline
\multicolumn{1}{|c|}{\multirow{8}{*}{\begin{tabular}[c]{@{}c@{}}\rotatebox[origin=c]{90}{Penetration Rate}\end{tabular}}} & \multirow{2}{*}{25\%}                       & \$1000                  & \multicolumn{1}{c|}{7242} & \multicolumn{1}{c|}{191} & 61   \\ \cline{3-6} 
\multicolumn{1}{|c|}{}                                                                            &                                             & \$10000               & \multicolumn{1}{c|}{7198} & \multicolumn{1}{c|}{14}  & 282  \\ \cline{2-6} 
\multicolumn{1}{|c|}{}                                                                            & \multirow{2}{*}{50\%}                       & \$1000                  & \multicolumn{1}{c|}{7463} & \multicolumn{1}{c|}{414} & 17   \\ \cline{3-6} 
\multicolumn{1}{|c|}{}                                                                            &                                             & \$10000                 & \multicolumn{1}{c|}{6975} & \multicolumn{1}{c|}{0}   & 519  \\ \cline{2-6} 
\multicolumn{1}{|c|}{}                                                                            & \multirow{2}{*}{75\%}                       & \$1000                 & \multicolumn{1}{c|}{6994} & \multicolumn{1}{c|}{500} & 0    \\ \cline{3-6} 
\multicolumn{1}{|c|}{}                                                                            &                                             & \$10000                 & \multicolumn{1}{c|}{6717} & \multicolumn{1}{c|}{0}   & 777  \\ \cline{2-6} 
\multicolumn{1}{|c|}{}                                                                            & \multicolumn{1}{l|}{\multirow{2}{*}{100\%}} & \$1000                  & \multicolumn{1}{c|}{6994} & \multicolumn{1}{c|}{500} & 0    \\ \cline{3-6} 
\multicolumn{1}{|c|}{}                                                                            & \multicolumn{1}{l|}{}                       & \$10000                 & \multicolumn{1}{c|}{6472} & \multicolumn{1}{c|}{28}  & 994  \\ \hline
\end{tabular}
\vspace{0.1cm}
\caption{Distribution of the offered incentives in Experiment~I with different penetration rates.} \label{table:distIncentiveI_7-8}

\end{table}

\vspace{-0.4cm}

\begin{table}[H]
\centering
\begin{tabular}{|c|c|c|c|}
\hline
\multirow{2}{*}{Budget} & \multicolumn{3}{c|}{Incentive} \\ \cline{2-4} 
                  &    \$0   &   \$2    &   \$10    \\ \hline
             \$1000     &   14645    & 435      &  13    \\ \hline
             \$10000     &   12614    &  1849     &  630     \\ \hline
\end{tabular}
\vspace{+0.15cm}
\caption{Distribution of the offered incentives in Experiment~II for incentive set $\cI_1$ with penetration rate of 100\%.}
\label{table:distIncentiveII_settingI_7-8}
\end{table}

\vspace{-0.4cm}

\begin{table}[H]
\centering
\begin{tabular}{|c|c|c|c|c|c|c|}
\hline
\multirow{2}{*}{Budget} & \multicolumn{6}{c|}{Incentive} \\ \cline{2-7} 
                  & \$0  & \$1  & \$2  & \$3  & \$5  & \$10  \\ \hline
               \$1000   & 14509  & 351  & 152  & 30  & 51  & 0  \\ \hline
                \$10000  &  12682 & 184  & 305  & 832  & 838  & 252  \\ \hline
\end{tabular}
\vspace{+0.1cm}
\caption{Distribution of the offered incentives in Experiment~II for incentive set $\cI_2$ with penetration rate of 100\%.} 
\label{table:distIncentiveII_settingII_7-8}
\end{table}

\vspace{-0.4cm}

\begin{table}[H]
\centering
\begin{tabular}{|c|c|c|c|}
\hline               \\ \hline
\multirow{2}{*}{Budget} & \multicolumn{3}{c|}{Incentive} \\ \cline{2-4} 
                  &    \$0   &   \$2    &   \$10    \\ \hline
             \$1000     &   7720    & 500      &  0    \\ \hline
             \$10000     &   6916    &  380     &  924     \\ \hline
\end{tabular}
\vspace{+0.15cm}
\caption{Distribution of the offered incentives in Experiment~III for model~\eqref{eq:LP2} with penetration rate of 100\%.} 
\label{table:distIncentiveIII_Linear_7-8}
\end{table}

\vspace{-0.4cm}

\begin{table}[H]
\centering
\begin{tabular}{cc|c|ccc|}
\cline{3-6}
\multicolumn{2}{c|}{\multirow{2}{*}{}}                                                                                                          & \multirow{2}{*}{Budget} & \multicolumn{3}{c|}{Incentive}                              \\ \cline{4-6} 
\multicolumn{2}{c|}{}                                                                                                                           &                         & \multicolumn{1}{c|}{\$0}  & \multicolumn{1}{c|}{\$2} & \$10 \\ \hline
\multicolumn{1}{|c|}{\multirow{8}{*}{\begin{tabular}[c]{@{}c@{}}\rotatebox[origin=c]{90}{Penetration Rate}\end{tabular}}} & \multirow{2}{*}{25\%}                       & \$1000                  & \multicolumn{1}{c|}{8040} & \multicolumn{1}{c|}{100} & 80   \\ \cline{3-6} 
\multicolumn{1}{|c|}{}                                                                            &                                             & \$10000                 & \multicolumn{1}{c|}{7866} & \multicolumn{1}{c|}{115} & 239  \\ \cline{2-6} 
\multicolumn{1}{|c|}{}                                                                            & \multirow{2}{*}{50\%}                       & \$1000                  & \multicolumn{1}{c|}{7894} & \multicolumn{1}{c|}{285} & 41   \\ \cline{3-6} 
\multicolumn{1}{|c|}{}                                                                            &                                             & \$10000                 & \multicolumn{1}{c|}{7153} & \multicolumn{1}{c|}{94}  & 973  \\ \cline{2-6} 
\multicolumn{1}{|c|}{}                                                                            & \multirow{2}{*}{75\%}                       & \$1000                  & \multicolumn{1}{c|}{7821} & \multicolumn{1}{c|}{374} & 25   \\ \cline{3-6} 
\multicolumn{1}{|c|}{}                                                                            &                                             & \$10000                 & \multicolumn{1}{c|}{6960} & \multicolumn{1}{c|}{335} & 925  \\ \cline{2-6} 
\multicolumn{1}{|c|}{}                                                                            & \multicolumn{1}{l|}{\multirow{2}{*}{100\%}} & \$1000                  & \multicolumn{1}{c|}{7780} & \multicolumn{1}{c|}{415} & 17   \\ \cline{3-6} 
\multicolumn{1}{|c|}{}                                                                            & \multicolumn{1}{l|}{}                       & \$10000                 & \multicolumn{1}{c|}{6891} & \multicolumn{1}{c|}{417} & 912  \\ \hline
\end{tabular}
\vspace{+0.15cm}
\caption{Distribution of the offered incentives in Experiment~III with different penetration rates for model~\eqref{eq:ObjFunctFinal}, Algorithm~\ref{alg:ADMM-ForOurProblem}.} 
\label{table:distIncentiveIII_ADMM_7-8}
\end{table}

\begin{table}[H]
\centering
\begin{tabular}{cc|c|ccc|}
\cline{3-6}
\multicolumn{2}{c|}{\multirow{2}{*}{}}                                                                                                          & \multirow{2}{*}{Budget} & \multicolumn{3}{c|}{Incentive}                              \\ \cline{4-6} 
\multicolumn{2}{c|}{}                                                                                                                           &                         & \multicolumn{1}{c|}{\$0}  & \multicolumn{1}{c|}{\$2} & \$10 \\ \hline
\multicolumn{1}{|c|}{\multirow{8}{*}{\begin{tabular}[c]{@{}c@{}}\rotatebox[origin=c]{90}{Penetration Rate}\end{tabular}}} & \multirow{2}{*}{25\%}                       & \$1000                  & \multicolumn{1}{c|}{7468} & \multicolumn{1}{c|}{110} & 78   \\ \cline{3-6} 
\multicolumn{1}{|c|}{}                                                                            &                                             & \$10000                 & \multicolumn{1}{c|}{7329} & \multicolumn{1}{c|}{9}   & 320  \\ \cline{2-6} 
\multicolumn{1}{|c|}{}                                                                            & \multirow{2}{*}{50\%}                       & \$1000                  & \multicolumn{1}{c|}{7980} & \multicolumn{1}{c|}{175} & 65   \\ \cline{3-6} 
\multicolumn{1}{|c|}{}                                                                            &                                             & \$10000                 & \multicolumn{1}{c|}{7576} & \multicolumn{1}{c|}{59}  & 585  \\ \cline{2-6} 
\multicolumn{1}{|c|}{}                                                                            & \multirow{2}{*}{75\%}                       & \$1000                  & \multicolumn{1}{c|}{7972} & \multicolumn{1}{c|}{185} & 63   \\ \cline{3-6} 
\multicolumn{1}{|c|}{}                                                                            &                                             & \$10000                 & \multicolumn{1}{c|}{7246} & \multicolumn{1}{c|}{78}  & 896  \\ \cline{2-6} 
\multicolumn{1}{|c|}{}                                                                            & \multicolumn{1}{l|}{\multirow{2}{*}{100\%}} & \$1000                  & \multicolumn{1}{c|}{7896} & \multicolumn{1}{c|}{280} & 44   \\ \cline{3-6} 
\multicolumn{1}{|c|}{}                                                                            & \multicolumn{1}{l|}{}                       & \$10000                 & \multicolumn{1}{c|}{7022} & \multicolumn{1}{c|}{248} & 950  \\ \hline
\end{tabular}
\vspace{0.1cm}
\caption{Distribution of the offered incentives in Experiment~III with different penetration rates for model~\eqref{eq:ObjFunctFinal}, Gurobi.} 
\label{table:distIncentiveIII_Gurobi_7-8}
\end{table}

\begin{table}[H]
\centering
\begin{tabular}{|c|cccc|}
\hline
\multirow{2}{*}{Budget} & \multicolumn{4}{c|}{Penetration Rate}                                                        \\ \cline{2-5} 
                        & \multicolumn{1}{c|}{25\%}  & \multicolumn{1}{c|}{50\%}  & \multicolumn{1}{c|}{75\%}  & 100\% \\ \hline
\$100                   & \multicolumn{1}{c|}{3.61}  & \multicolumn{1}{c|}{3.90}  & \multicolumn{1}{c|}{3.93}  & 4.91  \\ \hline
\$1000                  & \multicolumn{1}{c|}{16.27} & \multicolumn{1}{c|}{24.77} & \multicolumn{1}{c|}{31.52} & 31.24 \\ \hline
\$10000                 & \multicolumn{1}{c|}{20.41}      & \multicolumn{1}{c|}{41.05} & \multicolumn{1}{c|}{50.09} & 56.57 \\ \hline
\end{tabular}
\vspace{+0.15cm}
\caption{Effect of the penetration rate on travel time decrease (hour) in Experiment~I.} 
\label{table:PenEffectHourEx1}
\end{table}

\begin{table}[H]
\centering
\begin{tabular}{|c|cccc|}
\hline
\multirow{2}{*}{Budget} & \multicolumn{4}{c|}{Penetration Rate}                                                            \\ \cline{2-5} 
                        & \multicolumn{1}{c|}{25\%}   & \multicolumn{1}{c|}{50\%}   & \multicolumn{1}{c|}{75\%}   & 100\%  \\ \hline
\$100                   & \multicolumn{1}{c|}{0.47\%} & \multicolumn{1}{c|}{0.51\%} & \multicolumn{1}{c|}{0.51\%} & 0.64\% \\ \hline
\$1000                  & \multicolumn{1}{c|}{2.12\%} & \multicolumn{1}{c|}{3.23\%} & \multicolumn{1}{c|}{4.11\%} & 4.07\% \\ \hline
\$10000                 & \multicolumn{1}{c|}{2.66\%} & \multicolumn{1}{c|}{5.35\%} & \multicolumn{1}{c|}{6.63\%} & 7.37\% \\ \hline
\end{tabular}
\vspace{+0.15cm}
\caption{Effect of the penetration rate on percentage of travel time decrease in Experiment~I.} 
\label{table:PenEffectEx1}
\end{table}

\begin{table}[H]
\centering
\begin{tabular}{|c|cccc|}
\hline
\multirow{2}{*}{Budget} & \multicolumn{4}{c|}{Penetration Rate}                                                            \\ \cline{2-5} 
                        & \multicolumn{1}{c|}{25\%}   & \multicolumn{1}{c|}{50\%}   & \multicolumn{1}{c|}{75\%}   & 100\%  \\ \hline
\$100                   & \multicolumn{1}{c|}{39.41}  & \multicolumn{1}{c|}{43.78}  & \multicolumn{1}{c|}{28.91}  & 26.01  \\ \hline
\$1000                  & \multicolumn{1}{c|}{156.32} & \multicolumn{1}{c|}{215.23} & \multicolumn{1}{c|}{220.24} & 234.57 \\ \hline
\$10000                 & \multicolumn{1}{c|}{164.80} & \multicolumn{1}{c|}{255.77} & \multicolumn{1}{c|}{313.13} & 339.96 \\ \hline
\end{tabular}
\vspace{+0.15cm}
\caption{Effect of the penetration rate on travel time decrease (hour) in Experiment~III, model~\eqref{eq:ObjFunctFinal}, Algorithm~\ref{alg:ADMM-ForOurProblem}.} 
\label{table:PenEffectHourEx3ADMM}
\end{table}

\begin{table}[H]
\centering
\begin{tabular}{|c|cccc|}
\hline
\multirow{2}{*}{Budget} & \multicolumn{4}{c|}{Penetration Rate}                                                              \\ \cline{2-5} 
                        & \multicolumn{1}{c|}{25\%}   & \multicolumn{1}{c|}{50\%}   & \multicolumn{1}{c|}{75\%}    & 100\%   \\ \hline
\$100                   & \multicolumn{1}{c|}{1.31\%} & \multicolumn{1}{c|}{1.46\%} & \multicolumn{1}{c|}{0.96\%}  & 0.87\%  \\ \hline
\$1000                  & \multicolumn{1}{c|}{5.20\%} & \multicolumn{1}{c|}{7.16\%} & \multicolumn{1}{c|}{7.33\%}  & 7.81\% \\ \hline
\$10000                 & \multicolumn{1}{c|}{5.48\%} & \multicolumn{1}{c|}{8.51\%} & \multicolumn{1}{c|}{10.42\%} & 11.31\% \\ \hline
\end{tabular}
\vspace{+0.15cm}
\caption{Effect of the penetration rate on percentage of travel time decrease in Experiment~III, model~\eqref{eq:ObjFunctFinal}, Algorithm~\ref{alg:ADMM-ForOurProblem}.} 
\label{table:PenEffectEx3ADMM}
\end{table}

\begin{table}[H]
\centering
\begin{tabular}{|c|cccc|}
\hline
\multirow{2}{*}{Budget} & \multicolumn{4}{c|}{Penetration Rate}\\ \cline{2-5} 
                        & \multicolumn{1}{c|}{25\%}   & \multicolumn{1}{c|}{50\%}   & \multicolumn{1}{c|}{75\%}   & 100\%  \\ \hline
\$100                   & \multicolumn{1}{c|}{70.75}  & \multicolumn{1}{c|}{77.67}  & \multicolumn{1}{c|}{86.99}  & 88.51  \\ \hline
\$1000                  & \multicolumn{1}{c|}{201.04} & \multicolumn{1}{c|}{288.96} & \multicolumn{1}{c|}{362.90} & 402.51 \\ \hline
\$10000                 & \multicolumn{1}{c|}{215.95} & \multicolumn{1}{c|}{361.06} & \multicolumn{1}{c|}{473.70} & 536.59 \\ \hline
\end{tabular}
\vspace{+0.15cm}
\caption{Effect of the penetration rate on travel time decrease (hour) in Experiment~III, model~\eqref{eq:ObjFunctFinal}, Gurobi.} 
\label{table:PenEffectHourEx3Gurobi}
\end{table}

\begin{table}[H]
\centering
\begin{tabular}{|c|cccc|}
\hline
\multirow{2}{*}{Budget} & \multicolumn{4}{c|}{Penetration Rate}                                                               \\ \cline{2-5} 
                        & \multicolumn{1}{c|}{25\%}   & \multicolumn{1}{c|}{50\%}    & \multicolumn{1}{c|}{75\%}    & 100\%   \\ \hline
\$100                   & \multicolumn{1}{c|}{2.35\%} & \multicolumn{1}{c|}{2.58\%}  & \multicolumn{1}{c|}{2.89\%}  & 2.95\%  \\ \hline
\$1000                  & \multicolumn{1}{c|}{6.69\%} & \multicolumn{1}{c|}{9.62\%}  & \multicolumn{1}{c|}{12.08\%} & 13.39\% \\ \hline
\$10000                 & \multicolumn{1}{c|}{7.19\%} & \multicolumn{1}{c|}{12.01\%} & \multicolumn{1}{c|}{15.76\%} & 17.86\% \\ \hline
\end{tabular}
\vspace{+0.15cm}
\caption{Effect of the penetration rate on percentage of travel time decrease in Experiment~III, model~\eqref{eq:ObjFunctFinal}, Gurobi.} 
\label{table:PenEffectEx3Gurobi}
\end{table}

\end{document}